\documentclass[BCOR=12mm, DIV=default,twoside, pagesize, draft]{scrartcl} 
\usepackage{enumitem}

\addtokomafont{caption}{\small}
\setkomafont{captionlabel}{\sffamily\bfseries}

\usepackage{tikz}  
\usepackage{amsmath, amssymb, amsthm}
\usepackage{accents}
\usetikzlibrary{arrows,automata}   
\usepackage{todonotes}
\usepackage[utf8x]{inputenc}
\usepackage{libertine} 
\usepackage[T1]{fontenc}
\usepackage[english]{babel}
\usepackage{url}  %\usepackage{hyperref,faktor,csquotes,color}
\usepackage{faktor,csquotes,color}
\usepackage{pgfplots} 
\usepackage{stackengine}
\usepackage{breqn}
\usepackage{graphicx}
\usepackage{float} 
\usetikzlibrary{datavisualization.formats.functions}
\pgfplotsset{compat=1.12}
\usepgfplotslibrary{fillbetween}
\usetikzlibrary{patterns}
\usepackage{dsfont}

\usepackage{changes}

\definechangesauthor[name={Kristoffer}, color=blue]{kri}
\newcommand{\stkout}[1]{\ifmmode\text{\sout{\ensuremath{#1}}}\else\sout{#1}\fi}
\setdeletedmarkup{\stkout{#1}}

\usepackage{units,varioref}

\usepackage[shortcuts]{extdash}

\DeclareMathOperator{\Span}{span}
\DeclareMathOperator{\rank}{rank}

\usepackage{scrlayer-scrpage}
\KOMAoptions{draft=false}

\usepackage{mathtools}
\mathtoolsset{showonlyrefs,showmanualtags}
\numberwithin{equation}{section}

\allowdisplaybreaks 
\clearpairofpagestyles
\automark[section]{section}
\ohead{\headmark}
\ofoot[\pagemark]{\pagemark}

\newtheoremstyle{break}{\topsep}{\topsep}{\itshape}{}{\bfseries}{.}{\newline}{}
\newtheoremstyle{exampl}{\topsep}{\topsep}{\upshape}{}{\bfseries}{.}{\newline}{}

\theoremstyle{plain}% default
\newtheorem{theorem}{Theorem}[section]
\newtheorem{lemma}[theorem]{Lemma}  
\newtheorem{proposition}[theorem]{Proposition}  
\newtheorem{corollary}[theorem]{Corollary} 
 
\newtheorem{remark}[theorem]{Remark} 
\theoremstyle{definition}

\newtheorem{example}[theorem]{Example}

\def\F{\mathcal F}

\def\E{\mathbb E}
\def\R{\mathbb{R}}

\def\calL{\mathcal{L}}
\def\N{\mathbb{N}}

\def\P{\mathcal P}

\newcommand{\D}{\mathcal D}
\newcommand{\C}{\mathcal C}

\def\P{{\mathbb P}}

\newtheorem{assumption}[theorem]{Assumption}

\title{On hypoellipticity of degenerate operators in testing and detection problems}

\author{Erhan Bayraktar\footnote{Department of Mathematics, University of Michigan, United States. 
		E-mail address: erhan@umich.edu. This author is partially supported by the National Science Foundation under grant DMS-2602036 and by the Susan M. Smith Chair.} and Yuqiong Wang\footnote{Department of Mathematics, University of Michigan, United States.
		E-mail address: yuqw@umich.edu. }}

\date{\today}

\begin{document}
	\maketitle  
	
\begin{abstract} 
We study a class of degenerate diffusion generators arising in sequential testing and quickest detection problems with partial information. The observation process is driven by $k$ independent Brownian motions, while the hidden state takes $n+1$ values, with $k<n$.
After transforming to posterior likelihood coordinates, we analyze H\"ormander's condition both in the absence of state switching (testing) and in its presence (detection). We characterize hypoellipticity in the testing case and give explicit sufficient conditions in the detection case. We further study the stationary posterior operator, its parabolic extension, the joint observation-posterior operator, and its parabolic extension; their H\"ormander conditions need not coincide. We characterize the relationships among these operators under our main structural regimes and discuss their probabilistic consequences and the regularity of the associated optimal stopping problem. \end{abstract}
\section{Introduction}
\label{sec_intro}
We investigate hypoellipticity for a class of degenerate diffusion operators arising from filtering problems, particularly those in quickest detection and sequential hypothesis testing. We consider a $k$-dimensional Brownian motion $X$ with unknown drift. The drift of $X$ is determined by a hidden continuous-time Markov process $\{\theta_t\}_{t\geq 0}$ with state space $\{0,1, \dots, n\}$ and infinitesimal generator $Q$. The process $\theta$ is independent of the driving Brownian motion $W$, and $X$ satisfies
\begin{equation}
\label{x_dynamics}
    dX_t = \sum_{j=0}^n 1_{\theta_t=j} \lambda_jdt+dW_t, \quad X_0 = 0.
\end{equation}
Here, $\lambda_i\in \R^k$, $i\in \{0,1, \dots, n\}$, are known vectors. After filtering, the posterior is $n$-dimensional but is driven by only $k$ Brownian directions. When $k<n$, the generator is degenerate on the probability simplex. We study when Lie brackets transfer the observed noise into the missing posterior directions.\par
Many statistical applications based on \eqref{x_dynamics} require accurate and timely decisions, such as declaring a regime change or choosing the most likely state while accounting for observation costs. Reference \cite{S} studies two problems involving one-dimensional Brownian motion: $(i)$ a sequential hypothesis testing problem that seeks to identify the drift quickly while accounting for observation costs, in which case $\theta_t = \theta_0 \in \{0,1\}$ and $Q = 0$; and $(ii)$ a quickest disorder detection problem that seeks to identify the drift-change time quickly while penalizing false alarms. In the latter case, $\theta_t\in\{0,1\}$ and
\[Q = \begin{bmatrix}
   -\lambda&\lambda\\
   0&0
\end{bmatrix}.
\]
In other words, the Brownian motion acquires a constant drift at an exponentially distributed time with rate $\lambda$, and the drift never reverts. Both problems can be written as optimal stopping problems of the form
\begin{equation}
    \inf_{\tau} \E\left[g(\Pi_{\tau})+ \int_0^{\tau} h(\Pi_t)dt\right]
\end{equation}
where the state process $\Pi$ is the posterior probability $\Pi_t: = \P(\theta_t = 1\vert \F_t^X)$. Here $g$ is the penalty incurred upon making a decision, and $h$ is the running cost incurred before stopping. These one-dimensional problems have various extensions and analogous Poisson formulations because detecting a change in the drift of a Brownian motion parallels detecting a change in the intensity of a Poisson process; see, for example, \cite{Peskir_poisson, Peskir_poisson2}. Reference \cite{EV} treats cases in which the prior distribution is not Bernoulli and the hypotheses are composite. The authors of \cite{campbell} study a problem with a random, $\theta$-dependent time horizon. More closely related to the present paper, \cite{B_poisson} and \cite{zhitlukhin} study problems in which the natural sufficient statistics are two-dimensional but there is only one Brownian source. In a higher-dimensional setting, \cite{dayanik} studies a multisource detection problem with independent Poisson and Brownian observations that can be viewed as a generalized one-dimensional formulation. For a similar reduction of a multichannel problem to one dimension, see \cite{BR}. Reference \cite{BP_detection} studies a multidimensional Poisson problem with a unified stopping rule. References \cite{B_poisson} and \cite{EW2} study independent Poisson and Brownian sources, respectively, with change points that occur independently in each direction. More recently, \cite{gapeev} considers a multidimensional detection problem with correlated Brownian motions.\par
Beyond purely statistical applications, many problems feature payoffs that depend directly on the underlying process, so partial information about the unknown state affects the payoff only implicitly. This framework has a broad range of related applications. In such cases, the value function can be formulated in terms of the observation process and the belief process $\Pi$ obtained by filtering. See, for example, \cite{decamps} for an investment problem with incomplete information and \cite{gapeev2,guo} for the pricing of American options with regime switching. We refer to \cite{EW1} for a formulation in which the payoff may depend on time, the unknown state $\theta$, and the observation process, while the time horizon may be random and $\theta$-dependent. As discussed in Section \ref{sec_consequence}, our analysis extends to cases in which the payoff depends on both $\theta$ and $x$. We also refer to \cite{BL_inventory} for an impulse control problem. Although our results are formulated in terms of optimal stopping, they also have implications for singular and impulse control problems.\par
Our study of this multidimensional testing and detection problem in the degenerate case is motivated by \cite{caffarelli}, in which the authors study the problem of testing the state of $\theta$ with a specific payoff structure when $k=n$. In this setting, the operator is uniformly elliptic in the interior and degenerates only on the boundary. In Section $11$ of \cite{caffarelli}, the authors briefly consider the case $k<n$ with $Q=0$ and derive structural properties using probabilistic arguments. Inspired by this work, we focus on the case $k<n$, in which the operator is degenerate everywhere and standard theories requiring uniform ellipticity do not apply. Nevertheless, hypoellipticity may still hold if the missing diffusion directions can be recovered through Lie brackets. Thus, smoothness and other regularity properties may hold even when ellipticity fails. We identify settings in which this mechanism restores regularity by checking \emph{H\"ormander's finite-rank condition}.
\par 
Only a few works in the filtering literature discuss the hypoellipticity of the associated operators. In \cite{peskir_weak}, the authors show that when H\"ormander's condition holds for the differential operator of a strong Markov process, a weak solution can be upgraded to a smooth solution whenever the source term is smooth. In \cite{peskir_detection}, the authors solve a multidimensional disorder detection problem in which $K$ of $N$ underlying processes acquire a drift simultaneously, making the operator degenerate elliptic. We discuss this problem in Example \ref{eg2}.\par
 
Our main contribution is to identify conditions under which the degenerate generator arising in sequential testing and quickest detection remains hypoelliptic. We first give an exact characterization of the testing case. In Theorem \ref{thm_testing}, we show that when $Q=0$, H\"ormander's condition holds if and only if the drift-difference matrix has full rank and the vector of squared norms of the drift differences does not lie in its row space. In particular, hypoellipticity fails automatically when $n>k+1$. We then turn to the detection case, in which $Q\neq 0$. Transitions of the hidden state introduce additional terms into the posterior drift, and repeated commutation with the diffusion fields can generate directions that are unavailable in the testing problem. In Theorem \ref{thm_suff_1} and Proposition \ref{thm_suff_2}, we give two explicit, parameter-based criteria for determining whether these additional directions span the posterior space.\par
Our next contribution is to distinguish the H\"ormander conditions for the various operators arising from the same filtering model. The relevant operator depends on the stopping or control problem under consideration. We primarily consider the stationary operator in the posterior likelihood $\phi$-coordinate, denoted by $\calL_{\phi}$. Finite-horizon problems are governed by $-\partial_t+\calL_\phi$, while problems in which the payoff depends on the observation process are governed by the joint generator $\calL_{\phi,x}$. Finally, time-dependent problems whose payoffs depend on the observation process are governed by the parabolic joint operator $-\partial_t+\calL_{\phi,x}$. The H\"ormander conditions for these four operators are genuinely different. In particular, stationary hypoellipticity need not imply the parabolic condition, and the parabolic condition may hold even when the joint operator fails to satisfy H\"ormander's condition. We give sufficient conditions for the parabolic and joint operators to be hypoelliptic and identify when these conditions fail, both when $Q=0$ and when $Q\neq 0$. We also show that all four operators are hypoelliptic under the assumptions of Theorem \ref{thm_suff_1}.\par
We further discuss the probabilistic and optimal-stopping consequences of hypoellipticity. The parabolic H\"ormander condition implies the strong Feller property in the interior, while behavior on the boundary depends on the transition structure of $Q$. For the associated stopping problem, hypoellipticity yields smoothness of the value function in the continuation region under suitable regularity assumptions on the data but does not by itself imply smooth fit across the stopping boundary.\par
The rest of the paper is organized as follows. In Section \ref{sec_formulation}, we formulate the stopping problem in terms of the posterior probability process and introduce the necessary notation. In Section \ref{sec_change}, we transform to the posterior likelihood process and state the associated generator. In Section \ref{sec_hypoellipticity}, we study hypoellipticity in two cases. When the infinitesimal generator is $Q = 0$, which we call the ``testing case,'' we characterize when H\"ormander's condition holds and show that it necessarily fails when $n>k+1$. When $Q \neq 0$, which we call the ``detection case,'' we give two sufficient conditions for H\"ormander's condition, each expressed as a simple parametric check on the underlying system. In Section \ref{sec_consequence}, we discuss the hypoellipticity of the parabolic and joint operators, their implications for the underlying process, and the regularity of the solution to the stopping problem in the interior of the continuation region. In Section \ref{sec_example}, we conclude with several motivating examples.

\section{Problem formulation}
\label{sec_formulation}
Define $\bar n : = \{0,1,\dots, n\}$, $[n]: = \{1,\dots, n\}$, $\bar k : = \{0,1,\dots, k\}$, and $[k]: = \{1,\dots, k\}$, with $k\leq n$. We work on a probability space $(\Omega, \F, \P_{\pi})$ supporting
\begin{itemize}
    \item a Markov process $(\theta_t)_{t\geq 0}$ taking values in $\bar n$ with infinitesimal generator $Q=(q_{ij})_{i,j\in \bar n}$, where
    $q_{ij}\geq 0$ if $i\neq j$, $q_{ii}\leq 0$, and $\sum_{j=0}^n q_{ij} =0$ for all $i\in\bar n$,
    \item a $k$-dimensional standard Brownian motion $W = (W^1,\dots, W^k)$ independent of $\theta$.
\end{itemize}
Let $\pi = (\pi_0,\pi_1, \dots, \pi_n)$ be a prior probability distribution satisfying $\sum_{i=0}^n \pi_i = 1$, where $\pi_i$ is the prior probability of the initial state $\theta_0=i$. Then $\P_{\pi}  = \sum_{i=0}^n \pi_i \P_i$, where $\P_i$ is a probability measure under which $\theta_0 = i$, and $\P_{\pi}(\theta_0 = i) = \pi_i$ for $i\in \bar n$. We denote expectation with respect to $\P_{\pi}$ by $\E_{\pi}$. Let $\lambda = (\lambda_0, \lambda_1,\dots, \lambda_n)$, where $\lambda_i \in \R^k$ for $i\in \bar n$. We assume that $\theta$ is unobservable and that we observe only a $k$-dimensional continuous process $X = (X^1, \dots, X^k)$. The process $X$ is driven by the Brownian motion $W$, and its drift is determined by the state of $\theta$:
\[dX_t = \sum_{j=0}^n 1_{\theta_t=j} \lambda_jdt+dW_t, \quad X_0 = 0.\] 
Thus, $X$ is driven by a $k$-dimensional Brownian motion and has $n+1$ possible drifts, with drift $\lambda_j$ when $\theta_t=j$. We denote by $\F_t^X$ the $\sigma$-algebra generated by $(X_s)_{0\leq s\leq t}$ and write $\F_t: = \F_t^X$. In the statistical interpretation, hypothesis $H_i$ is true at time $t$ if $\theta_t = i$. When the infinitesimal generator $Q$ is the zero matrix, $\theta_t = \theta_0$ for all $t\geq 0$. We refer to this as the \emph{sequential testing} problem and discuss it in detail in Section \ref{sec_hypoellipticity}.\par
Define the posterior probability process $\Pi = (\Pi^0,\dots, \Pi^n)$ by 
\[\Pi_t^i: = \P_{\pi} (\theta_t = i\vert \F_t), \quad t\geq 0, \enskip i\in \bar n.\]
The process $\Pi$ takes values in the $n$-dimensional simplex
\[P_{n+1} = \{\pi= (\pi_0, \pi_1,\dots, \pi_n)\in \R^{n+1}: \pi_i\geq 0, \sum_{i=0}^n \pi_i=1\}.\]
Thus $\Pi_t^i$ is the conditional probability that the hidden state equals $i$ given the observations up to time $t$. By standard filtering theory (see, e.g., \cite{LS}), for each $j\in \bar n$, the posterior coordinate $\Pi^j$ satisfies
\begin{equation}\label{filtering_pi}
    d\Pi_t^j = \sum_{i=0}^n q_{ij} \Pi_t^i dt+ \Pi_t^j (\lambda_j - \bar\lambda_t)\cdot  d\tilde W_t, \quad \Pi_0^j = \pi_j.
\end{equation}
Here $\bar\lambda_t := \sum_{i=0}^n \lambda_i\Pi_t^i \in \R^k$, and $\tilde W$ is the innovation process with $\tilde W_t = X_t -\int_0^t \bar \lambda_s ds$. Intuitively, the drift term describes the Markov transitions of $\theta$, while the martingale term updates the posterior using the innovation, that is, the new information in the observations. The process $\Pi$ is a strong Markov process with continuous paths, and its infinitesimal generator is given by
\begin{equation}
\label{pi-generator}
    \calL_{\pi} = \frac{1}{2} \sum_{i,j=0}^n \pi_i\pi_j (\lambda_i-\sum_{k=0}^n \pi_k\lambda_k)\cdot (\lambda_j-\sum_{l=0}^n \pi_l\lambda_l)\frac{\partial^2 }{\partial \pi_i \partial \pi_j}+\sum_{i,j=0}^n q_{ij} \pi_i\frac{\partial }{\partial \pi_j}.
\end{equation}
Let $g: [0,\infty)\times P_{n+1}\to \R$ be the immediate payoff upon stopping, $h: [0,\infty)\times P_{n+1}\to \R$ the running payoff accumulated before stopping, and $r: P_{n+1}\to [0,\infty)$ the discount rate. We assume that $g$, $h$, and $r$ are continuous, that $g$ is bounded, and that the following integrability condition holds:
 \[\E_{\pi}\left[ \int_0^{\infty}e^{-\int_0^{t}r(\Pi_s)ds}\left\vert  h(\Pi_t) \right\vert dt\right]<\infty. \]
Let $\mathcal{T}$ be the set of all $\F$-stopping times. We consider the following stopping problem:
\begin{equation}
\label{value_fcn}
    V = \sup_{\tau\in \mathcal{T}} \E_{\pi}\left[e^{-\int_0^{\tau}r(\Pi_s)ds}g(\Pi_{\tau}) + \int_{0}^{\tau} e^{-\int_0^{t}r(\Pi_s)ds}h(\Pi_t) dt\right].
\end{equation}
In \cite{caffarelli}, the authors study the case $k=n$, assuming that $\lambda_1-\lambda_0, \lambda_2-\lambda_0,\dots, \lambda_n-\lambda_0$ are linearly independent. The operator is then nondegenerate throughout the interior of the simplex $P_{n+1}$ and degenerates only on its boundary. If $k<n$, however, the operator is degenerate elliptic throughout the interior. We focus on this degenerate case and verify the hypoellipticity of the resulting operator by checking H\"ormander's condition \cite{homander}.\par
To check H\"ormander's condition, we write the infinitesimal generator $\calL_{\pi}$ as $\calL_{\pi} = D_0+\frac{1}{2}\sum_{r = 1}^k D_r^2 $, where $D = \{D_0,D_1, \dots, D_k\}$ consists of vector fields in $\R^n$ with $C^{\infty}$ coefficients. Recall that the Lie bracket of two fields $D_i$ and $D_j$ is defined by $[D_i, D_j] = D_iD_j - D_jD_i$. For an arbitrary multi-index $\alpha = (\alpha_1,\dots, \alpha_l)$, where $\alpha_i \in\bar k $ and $\vert\alpha\vert = l$, we define the iterated Lie brackets
\[D_{\alpha}: =[D_{\alpha_l},[D_{\alpha_{l-1}},\dots,[D_{\alpha_2}, D_{\alpha_1}]]].\]
We say that the system $D$ satisfies H\"ormander's finite-rank condition of order $s$ if
\[\operatorname{Lie} (D_0, \dots, D_k) := \{D_{\alpha}: \vert\alpha\vert \leq s\}\]
spans $\R^n$ at every point. It then follows that the operator $\calL_{\pi}$ is hypoelliptic \cite{homander}.
\begin{assumption}
    Throughout the paper, we assume the following:
    \begin{enumerate}
        \item $k<n$.
        \item The vectors $\lambda_0,\lambda_1,\dots,\lambda_n$ are pairwise distinct, and $\lambda_0$ is a vertex of their convex hull.
    \end{enumerate}
\end{assumption}
The second assumption implies $a_i\neq 0$, $a_i\neq a_j$, and $a_i\neq -a_j$ for all $i\neq j$. Without loss of generality, we may label one vertex of the convex hull as $\lambda_0$. The assumption that all possible drifts are pairwise distinct is natural, since otherwise some states would be indistinguishable and the effective number of states would decrease. A direct calculation shows that the infinitesimal generator \eqref{pi-generator} can be written as
\begin{equation}
\label{sumofsquares}
    \calL_{\pi} = D_0 + \frac{1}{2}\sum_{i=1}^k D_i^2, 
\end{equation}
where $D_1,\dots,D_k$ are vector fields defined on the interior of $P_{n+1}$. In particular, $D_r = \sum_{i=0}^n \pi_i (\lambda_{ir} - \bar \lambda_{ir}(\pi)) \partial_{\pi_i}$ for $r\in[k]$, and $D_0 = C(\pi)\cdot \nabla_{\pi}$, where
\[C_j(\pi) = \sum_{i=0}^n q_{ij} \pi_i -\frac{1}{2} \sum_{r=1}^k \sum_{i=0}^n \pi_i (\lambda_{ir} - \bar \lambda_{ir}(\pi)) \partial_{\pi_i}\left(\pi_j (\lambda_{jr} - \bar \lambda_{jr}(\pi))\right)\]
for $j\in \bar n$. Algebraically, the Lie bracket of any two vector fields $D_i$ and $D_j$ generates nonaffine, globally coupled coefficients involving $\bar \lambda$. In the next section, we perform a change of coordinates and check H\"ormander's condition in the transformed coordinates.
\section{A change of coordinates} 
\label{sec_change}
For $i\in \bar n$, define the posterior likelihood ratio with respect to state $0$ by $\Phi^i_t : = \frac{\Pi_t^i}{\Pi_t^0}$ whenever $\Pi_t^0>0$, and set $\phi =(\phi^1,\dots, \phi_n)$ with $\phi_i = \frac{\pi_i}{\pi_0}$. Since $\Phi_t^0 \equiv 1$, define $Y_t = \sum_{i=0}^n\Phi_t^i=1+ \sum_{i=1}^n\Phi_t^i$. The posterior probability process then satisfies $\Pi_t^i = \frac{\Phi_t^i}{Y_t}$ for $i\in \bar{n}$. By standard filtering theory, the posterior coordinates satisfy \eqref{filtering_pi}. Applying It\^o's formula to $\Phi^i_t$ yields
\begin{align*}
    d\Phi^i_t = &\frac{1}{\Pi_t^0} d\Pi_t^i-\frac{\Phi_t^i}{\Pi_t^0} d\Pi_t^0+ \frac{\Phi_t^i}{(\Pi_t^0)^2}d\langle \Pi^0\rangle_t -\frac{1}{(\Pi_t^0)^2}d\langle \Pi^0,\Pi^i\rangle_t\\
     =& \sum_{m=0}^n q_{mi}\Phi_t^m dt + \Phi_t^i(\lambda_i-\bar \lambda_t)\cdot d\tilde W_t -\Phi_t^i\left(\sum_{m=0}^nq_{m0}\Phi_t^m dt +(\lambda_0-\bar \lambda_t)\cdot d\tilde W_t\right)\\
     &+(\Phi_t^i)  \|\lambda_0-\bar \lambda_t\|^2-\Phi_t^i (\lambda_i-\bar \lambda_t)\cdot(\lambda_0-\bar \lambda_t)dt\\
     =&  \left(\sum_{m=0}^n q_{mi}\Phi_t^m - \Phi_t^i \sum_{m=0}^n q_{m0}\Phi_t^m \right)dt- \Phi_t^i(\lambda_i-\lambda_0)\cdot(\lambda_0-\bar \lambda_t)dt +\Phi_t^i (\lambda_i-\lambda_0)\cdot d \tilde W_t.
     \end{align*}

The average drift is $ \bar \lambda_t = \frac{1}{Y_t}\left(\lambda_0+ \sum_{i=1}^n\lambda_i \Phi_t^i \right)$. Therefore, the term $(\lambda_i-\lambda_0)\cdot(\lambda_0-\bar \lambda_t)$ in the drift coefficient can be written as 
\begin{align*}
   (\lambda_i-\lambda_0)\cdot(\lambda_0-\bar \lambda_t)& = (\lambda_i-\lambda_0)\cdot\left(\lambda_0-\frac{1}{Y_t}\left(\lambda_0+ \sum_{m=1}^n\lambda_m \Phi_t^m \right)\right)\\
   & = (\lambda_i-\lambda_0)\cdot \left(\frac{(Y_t-1)\lambda_0- \sum_{m=1}^n\lambda_m \Phi_t^m}{Y}\right)\\
   & = (\lambda_i-\lambda_0)\cdot \left(\frac{1}{Y_t}\left( \sum_{m=1}^n(\lambda_0-\lambda_m )\Phi_t^m\right)\right)\\
   & = -\frac{1}{Y_t} \sum_{m=1}^n\Phi_t^m(\lambda_i-\lambda_0)\cdot (\lambda_m-\lambda_0).
\end{align*}
Substituting this expression into the dynamics of $\Phi_t^i$ gives
 \[d\Phi^i_t  =  \left(\sum_{m=0}^n \Phi_t^m (q_{mi}-q_{m0}\Phi_t^i)  + \frac{1}{Y_t} \sum_{m=1}^n\Phi_t^i \Phi_t^m(\lambda_i-\lambda_0)\cdot (\lambda_m-\lambda_0)\right)dt+\Phi_t^i (\lambda_i-\lambda_0)\cdot d \tilde W_t.\]
For $i\in \bar n$, define $a_i: = \lambda_i - \lambda_0\in \R^k$ to be the difference between the drifts in states $i$ and $0$. By our choice of $\lambda_0$ and the assumption that all drifts are pairwise distinct,
\[a_i\neq a_j,\enskip a_i \neq -a_j, \enskip \text{ for all distinct }i,j \in [n].\]
Let $A:=(a_1,\dots, a_n) \in \R^{k\times n}$ be the matrix whose columns are the drift differences, and set $\Sigma : = A^TA\in \R^{n\times n}$. Then $\Sigma_{ij}:=\Sigma_{i,j} = a_i\cdot a_j$, so $\Sigma$ is positive semidefinite and $\operatorname{rank } \Sigma\leq k$. The dynamics of the $i$th posterior likelihood process can then be written as 
\begin{equation}
\label{phi-sde}
   d\Phi^i_t  =  \left(\sum_{m=0}^n \Phi_t^m (q_{mi}-q_{m0}\Phi_t^i)  + \frac{1}{Y} \sum_{m=1}^n\Sigma_{i,m}\Phi_t^i \Phi_t^m\right)dt+\Phi_t^i a_i\cdot d \tilde W_t
\end{equation}   
with $\Phi^i_0 = \phi_i$. The infinitesimal generator $\calL$ of the process $\Phi$ is
\begin{equation}
\label{phi_operator}
 \calL = \frac{1}{2}\sum_{i,j=1}^n \Sigma_{ij} \phi_i\phi_j \frac{\partial^2}{\partial_{\phi_i}\partial_{\phi_j}} +\frac{1}{y(\phi)}\sum_{i,j=1}^n\Sigma_{ij} \phi_i\phi_j \frac{\partial}{\partial_{\phi_j}} + \sum_{j=1}^n\sum_{i=0}^n (q_{ij} - q_{i0}\phi_j)\phi_i\frac{\partial}{\partial_{\phi_j}},
\end{equation}
consistent with \cite{caffarelli}. The diffusion matrix in \eqref{phi_operator} has rank at most $k$ and is therefore generally degenerate. Below, we discuss the hypoellipticity of \eqref{phi_operator}. Denote the interior of $P_{n+1}$ by $\mathring{P}_{n+1}$. The coordinate transformation from $\mathring{P}_{n+1}$ to $(0, \infty)^n$ is a $C^{\infty}$ diffeomorphism. Therefore, hypoellipticity in the $\phi$-coordinate on $(0,\infty)^n$ implies hypoellipticity in the $\pi$-coordinate on $\mathring{P}_{n+1}$. Given the initial value $\Phi_0 = (\phi_1,\dots, \phi_n) \in (0,\infty)^n$, we can find an explicit solution to the linear SDE \eqref{phi-sde}. Define
\[b_1^i(\phi) = q_{ii} - \sum_{m=0}^n q_{m0}\phi_m +\frac{1}{y(\phi)}\sum_{m=1}^n \Sigma_{im} \phi_m,\quad b_2^i(\phi) = \sum_{m\neq i}^n q_{mi}\phi_m, \]
and define the process $Z_t^i: = \exp\left\{\int_0^t b_1^i(\Phi_s)ds +a_i\cdot \tilde{W}_t - \frac{1}{2} \|a_i\|^2 t\right\}$. Then, for $i\in [n]$, $\Phi_t^i$ is given by
\[\Phi_t^i = Z_t^i\left(\phi_i + \int_0^t (Z_s^i)^{-1} b_2^i(\Phi_s)ds\right).\]
Since $Z_t^i>0$ for all $i$, it follows that $\Phi_t\in(0,\infty)^n$ for all $t\geq 0$. Equivalently, if $\pi \in \mathring{P}_{n+1}$, then $\Pi_t \in\mathring{P}_{n+1} $ for all $t\geq 0$. Now suppose that $\pi \in \partial {P}_{n+1}$ and $\pi_i = 0$ for some $i$. If $b_2^i(\phi) = \sum_{m\neq i}^n q_{mi}\phi_m>0$, then $\Phi_t^i>0$ for all $t>0$. If this condition holds for every coordinate that is initially zero, then the posterior process belongs to $\mathring P_{n+1}$ for every $t>0$. Even if $b_2^i(\phi)=0$ at the initial point, the coordinate $\Phi^i$ may become positive through a sequence of transitions governed by $Q$. In the testing case, however, $Q=0$ and $b_2^i\equiv 0$ for all $i$, so a process starting on the boundary remains there. If $r$ prior probabilities are positive, the stopping problem reduces to an $(r-1)$-dimensional problem.
\section{Hypoellipticity of $\calL_{\phi}$}
\label{sec_hypoellipticity}
In this section, we examine H\"ormander's condition for the operator \eqref{phi_operator} and give conditions under which the operator is hypoelliptic. Specifically, we determine whether iterated Lie brackets involving the $k$ diffusion directions, the drift, and, when present, the state-transition term associated with $Q$ span all $n$ directions. We begin with the special case in which $\theta$ does not depend on $t$, or equivalently, $Q = 0$. The observed process has a constant drift determined by the $(n+1)$-valued random variable $\theta$. We say that hypothesis $H_j$ is true if $\theta = j$; in statistical applications, one may wish to test hypotheses $H_0,\dots, H_n$. We call this the \emph{testing case} and give an exact classification.\par
We then consider the general case $Q\neq 0$. As mentioned in Section \ref{sec_intro}, the classical one-dimensional disorder detection problem assumes that $$Q = \begin{bmatrix}
    -\lambda &\lambda\\
    0&0
\end{bmatrix},$$ so a drift change is permanent once it occurs. We allow a general Markov generator matrix $Q$, so the drift may return to state $0$ or jump to other states at exponentially distributed times. In some applications, one may wish to detect the exact times of these changes. We therefore call this the \emph{detection case}, although it is more general than the classical detection problem. We give sufficient conditions for hypoellipticity and discuss their necessity.\par

\subsection{The testing case}
In the absence of state switching, the process $\Phi_t^i$ has the following dynamics: 
\[d\Phi_t^i= \Phi_t^i\sum_{l=1}^k a_{il}dW_t^l +\frac{1}{Y_t} \sum_{j=1}^n \Sigma_{ij}\Phi_t^i \Phi_t^j dt, \]
where the process $Y$ is defined above and $a_{il}$ is the $l$th component of the vector $a_i: = \lambda_i-\lambda_0\in \R^k$. Consistent with \cite{caffarelli}, we can write the operator in the $\phi$-coordinate as \eqref{phi_operator}. We first identify suitable vector fields that put the generator $\calL$ into the ``sum of squares'' form in \eqref{sumofsquares}.
Introduce the $k$ first-order fields $D_r = \sum_{i=1}^n a_{ir} \phi_i \partial_{\phi_i}$ for $r \in[k]$. The field $D_r$ is the diffusion direction associated with the $r$th Brownian coordinate, weighted by the current likelihood ratios. We compute
\begin{align*}
    D_r^2 & = \sum_{i=1}^n a_{ir} \phi_i\partial_{\phi_i}\left(\sum_{j=1}^n a_{jr} \phi_j\partial_{\phi_j}\right)\\
    & = \sum_{i=1}^n a_{ir} \phi_i\left(\sum_{j\neq i}^n a_{jr} \phi_j\partial^2_{\phi_i\phi_j}+a_{ir}\left(\phi_i\partial^2_{\phi_i\phi_i}+\partial_{\phi_i}\right)\right)\\
    & = \sum_{i=1}^n a^2_{ir} \phi_i\partial_{\phi_i} +\sum_{i,j=1}^n a_{ir}a_{jr}\phi_i\phi_j\partial^2_{\phi_i\phi_j}.
\end{align*}
Summing over $r$, we obtain 
\begin{align*}
    \sum_{r=1}^k  D_r^2& = \sum_{r=1}^k\left(\sum_{i=1}^n a^2_{ir} \phi_i\partial_{\phi_i} +\sum_{i,j=1}^n a_{ir}a_{jr}\phi_i\phi_j\partial^2_{\phi_i\phi_j}\right)\\
    & = \sum_{i=1}^n\left(\sum_{r=1}^k a_{ir}^2\right)\phi_i\partial_{\phi_i} +\sum_{i,j=1}^n \left(\sum_{r=1}^k a_{ir}a_{jr}\right)\phi_i\phi_j\partial^2_{\phi_i\phi_j}\\
    & = \sum_{i=1}^n||a_i||^2\phi_i\partial_{\phi_i} +\sum_{i,j=1}^n \Sigma_{ij}\phi_i\phi_j\partial^2_{\phi_i\phi_j}.
\end{align*}
Define $D_0 = \frac{1}{y(\phi)}\sum_{i,j=1}^n \Sigma_{ij}\phi_i\phi_j \partial_{\phi_i} - \frac{1}{2}\sum_{i=1}^n||a_i||^2\phi_i\partial_{\phi_i}$. We can write the operator $\calL$ as
\[\calL = D_0+\frac{1}{2}\sum_{r=1}^k D_r^2.\]
Before proceeding, consider first-order vector fields $F$ and $G$ of the form $F = \sum_{i=1}^n f_i(\phi) \partial_{\phi_i}$ and $G = \sum_{i=1}^n g_i(\phi)\partial_{\phi_i}$, where $f_i, g_i: \R^n\to \R$ are $C^{\infty}$ functions. Their Lie bracket is
\begin{align*}
    [F,G] &=\sum_{i=1}^n f_i(\phi)\partial_{\phi_i}\left(\sum_{j=1}^n g_j(\phi) \partial_{\phi_j}\right) - \sum_{i=1}^n g_i(\phi)\partial_{\phi_i}\left(\sum_{j=1}^n f_j(\phi) \partial_{\phi_j}\right)\\
     = &\sum_{i,j=1}^n f_i(\phi)(\partial_{\phi_i} g_j)(\phi) \partial_{\phi_j}+\sum_{i,j=1}^n f_i(\phi) g_j(\phi) \partial_{\phi_i\phi_j}\\
    &-\sum_{i,j=1}^n g_i(\phi)(\partial_{\phi_i} f_j)(\phi) \partial_{\phi_j}-\sum_{i,j=1}^n g_i(\phi) f_j(\phi) \partial_{\phi_i\phi_j}\\
    =&\sum_{i=1}^n\left(\sum_{j=1}^n f_j(\phi)(\partial_{\phi_j} g_i)(\phi)-g_j(\phi)(\partial_{\phi_j} f_i)(\phi) \right)\partial_{\phi_i}\\
    & = \sum_{i=1}^n\left(F\left(g_i\right) -G\left(f_i\right) \right)\partial_{\phi_i}.
\end{align*}
\begin{lemma}
    \label{testing_bracket}
    Assume $Q=0$. Then the diffusion fields commute: $[D_r, D_u] = 0$ for all $r,u\in [k]$. In addition, for any $r\in [k]$, $[D_0, D_r]$ lies in the span of the diffusion fields:
    \[[D_0, D_r] = -\sum_{s=1}^k \left(\frac{\alpha_{rs}}{y(\phi)} - \frac{\beta_r\beta_s}{y^2(\phi)}\right) D_s,\]
    where $\alpha_{rs} = \sum_{j=1}^n a_{jr} a_{js} \phi_j$, $\beta_r =\sum_{j=1}^n a_{jr}\phi_j$, and $\beta_s=\sum_{j=1}^n a_{js}\phi_j$.
\end{lemma}
\begin{proof}
    First, observe that for any $r,u \in[k]$, the Lie bracket
$[D_r, D_u] = D_rD_u-D_uD_r = \sum_{i=1}^n \left(D_r(a_{iu}\phi_i) - D_u(a_{ir}\phi_i) \right)$. Indeed, $D_r(a_{iu}\phi_i) = a_{iu}D_r(\phi_i) = a_{iu}\sum_{j=1}^n a_{jr}\phi_j\partial_{\phi_j}(\phi_i) =a_{iu}\sum_{j=1}^n a_{jr}\phi_j1_{i=j} =  a_{iu}a_{ir}\phi_i$. Similarly, $D_u(a_{ir}\phi_i) = a_{ir}a_{iu}\phi_i$. Therefore, $D_r(a_{iu}\phi_i)-D_u(a_{ir}\phi_i)  = 0$ for every $i$, and $[D_r, D_u] = 0$. Thus, all the diffusion vector fields commute, and all iterated Lie brackets involving only diffusion fields vanish.\par
Similarly, fix any $r\in [k]$ and compute the bracket $[D_0, D_r]$. Write $D_0 = \sum_{i=1}^n f_i(\phi) \partial_{\phi_i}$ and $D_r = \sum_{i=1}^n g_i(\phi) \partial_{\phi_i}$. Then $f_i(\phi) = \phi_i \left(\frac{1}{y(\phi)} \sum_{l=1}^n \Sigma_{li}\phi_l -\frac{1}{2} ||a_i||^2\right)$
and $g_i(\phi) = a_{ir}\phi_i$. First note that $D_0(g_i(\phi)) = a_{ir}f_i(\phi)$. For $j\neq i$, taking the $j$th derivative of $f_i(\phi)$, we have
\[\partial_{\phi_j}(f_i(\phi))= \phi_i \partial_{\phi_j} \left(\frac{1}{y(\phi)} \sum_{j=1}^n \Sigma_{ij}\phi_j -\frac{1}{2} ||a_i||^2\right) = \phi_i  \left(\frac{\Sigma_{ij} y(\phi) - \sum_{u=1}^n \Sigma_{iu}\phi_u }{y^2(\phi)}\right).\]
For the diagonal term,
\[\partial_{\phi_i}(f_i(\phi)) = \phi_i \left(\frac{\Sigma_{ii} y(\phi) - \sum_{u=1}^n \Sigma_{iu}\phi_u }{y^2(\phi)}\right) +\frac{f_i(\phi)}{\phi_i}\]
which yields
\[D_r(f_i(\phi)) = a_{ir}f_i(\phi)+\phi_i\sum_{j=1}^n a_{jr}\phi_j\left(\frac{\Sigma_{ij} y(\phi) - \sum_{u=1}^n \Sigma_{iu}\phi_u }{y^2(\phi)}\right).\]
Combining $D_0(g_i(\phi))$ and $D_r(f_i(\phi))$, we find that the $i$th component of $[D_0, D_r]$ is
\[[D_0, D_r]^i  = -\phi_i\sum_{j=1}^n a_{jr}\phi_j\left(\frac{\Sigma_{ij} y(\phi) - \sum_{u=1}^n \Sigma_{i,u}\phi_u }{y^2(\phi)}\right) \partial_{\phi_i}. \] 
Recall that $\Sigma_{ij} = \sum_{s=1}^k a_{is}a_{js}$. Substituting this identity into the previous equation gives
\begin{align*}
    [D_0, D_r]^i & =  -\phi_i\sum_{j=1}^n a_{jr}\phi_j\left(\frac{\sum_{s=1}^k a_{is}a_{js}y(\phi) - \sum_{u=1}^n \sum_{s=1}^k a_{is}a_{us}\phi_u }{y^2(\phi)}\right) \partial_{\phi_i}\\
    & = -\sum_{s=1}^k\phi_i \left( \sum_{j=1}^n \frac{a_{jr}\phi_j}{y^2(\phi)} (y(\phi)a_{is}a_{js} -\sum_{u=1}^n a_{is}a_{us} \phi_u )\right)\partial_{\phi_i}\\
    & = -\sum_{s=1}^k a_{is}\phi_i \left(\frac{\alpha_{rs}}{y(\phi)} - \frac{\beta_r\beta_s}{y^2(\phi)}\right),
\end{align*}
where $\alpha_{rs} = \sum_{j=1}^n a_{jr} a_{js} \phi_j,\quad \beta_r =\sum_{j=1}^n a_{jr}\phi_j, \quad  \beta_s=\sum_{j=1}^n a_{js}\phi_j$ are smooth functions of $\phi$. Therefore, the Lie bracket 
\[ [D_0, D_r]  = \sum_{i=1}^n [D_0, D_r]^i = -\sum_{i=1}^n \sum_{s=1}^k a_{is}\phi_i \left(\frac{\alpha_{rs}}{y(\phi)} - \frac{\beta_r\beta_s}{y^2(\phi)}\right) = -\sum_{s=1}^k \left(\frac{\alpha_{rs}}{y(\phi)} - \frac{\beta_r\beta_s}{y^2(\phi)}\right) D_s.\]
Thus, $[D_0, D_r]\in \operatorname{span}_{C^{\infty}}\{D_1,\dots, D_k\} $ for any $r\in[k]$.
\end{proof}
\begin{remark}
    To interpret the factors $\alpha_{rs}$, $\beta_r$, and $\beta_s$, view $a_{ir}$ as the effect of the $r$th Brownian coordinate on the $i$th hypothesis. Recall that $\pi_i = \frac{\phi_i}{y(\phi)}$. The coefficient of the bracket $[D_0, D_r]$ in the $s$ direction can then be written as 
    \[\frac{\alpha_{rs}}{y(\phi)} - \frac{\beta_r\beta_s}{y^2(\phi)} = \sum_{j=1}^n a_{jr} a_{js} \pi_j - \sum_{j=1}^n a_{jr} \pi_j  \sum_{j=1}^na_{js} \pi_j,\]
    which is precisely the posterior covariance of the effects of the Brownian motion in the $r$ and $s$ directions. In particular, the bracket emphasizes directions with large covariance in absolute value and points in the opposite direction.
\end{remark}

\begin{lemma}
\label{testing_bracket_2}
Define $\mathcal{M}: =\operatorname{span}_{C^{\infty}}(D_0, D_1,\dots, D_k) = \left\{\sum_{r=0}^k f_r(\phi) D_r: f_r \in C^{\infty}((0,\infty)^n),\text{ for }r\in \bar k\right\}$. Then $Lie(D_0,D_1,\dots, D_k ) \subset \mathcal{M}$. Moreover, $\operatorname{dim } Lie(D_0,D_1,\dots, D_k)\leq k+1$.
\end{lemma}
\begin{proof}
The drift field and all diffusion fields belong to $\mathcal{M}$. Moreover, $[D_0, D_0] = [D_r, D_u] = 0\in \mathcal{M}$, and $[D_0, D_u] \in \mathcal{M}$ by Lemma \ref{testing_bracket}. Take any $M = \sum_{r=0}^k f_r(\phi) D_r \in \mathcal{M}$. We verify that $[D_0, M] \in \mathcal{M}$ and $[D_u, M] \in \mathcal{M}$ for any $u\in[k]$. For a fixed $u$, the Leibniz rule gives $[D_u, M] = \sum_{r = 1}^k [D_u, f_r(\phi) D_r]
 = \sum_{r = 1}^k D_u \left(f_r(\phi)\right) D_r -\sum_{r = 1}^k f_r(\phi) [D_r, D_u]  = \sum_{r = 1}^k D_u \left(f_r(\phi)\right) D_r \in \mathcal{M}$, by Lemma \ref{testing_bracket}. For the drift part, let us first write $\left[D_0, D_r\right]=\sum_{s=1}^k c_s^r(\phi) D_s$ for simplicity. Similarly,
\begin{align*}
[D_0, M] & = \sum_{r = 1}^k D_0 \left(f_r(\phi)\right) D_r -\sum_{r = 1}^k f_r(\phi) [D_r, D_0] \\
& = \sum_{r = 1}^k \left(D_0 \left(f_r(\phi)\right) D_r + f_r(\phi)\sum_{s=1}^k c_s^r(\phi) \right)D_r\in \mathcal{M}.
\end{align*}

Thus, finite iterated Lie brackets involving $D_0$ and the diffusion fields $D_1, \dots, D_k$ create no directions outside the span of $\{D_0,\dots, D_k\}$. Consequently, $\mathcal M$ is closed under these bracketing operations and is therefore a Lie subalgebra. Hence $Lie(D_0,D_1,\dots, D_k ) \subset \mathcal{M}$, and $\operatorname{dim } Lie(D_0,D_1,\dots, D_k)\leq k+1$.
\end{proof}
This analysis shows that, \emph{in the testing case}, H\"ormander's condition may hold when $n = k+1$. When $n>k+1$, however, hypoellipticity is impossible. The following theorem characterizes precisely when the operator is hypoelliptic.
\begin{theorem}
\label{thm_testing}
   Assume $Q=0$ and define $A: = (a_1, a_2,\dots, a_n)\in \R^{k\times n}$. For $n=k+1$, H\"ormander's condition holds if and only if $\operatorname{rank}\left(A\right) = k$ and the vector $\left(||a_1||^2, \dots, ||a_n||^2\right)$ does not belong to $\operatorname{rowsp}(A)$. If $n>k+1$, H\"ormander's condition fails.
\end{theorem}
\begin{proof}
    The failure of H\"ormander's condition when $n>k+1$ follows immediately from Lemma \ref{testing_bracket_2}. Now suppose that $n=k+1$, and let $A_r$ denote the $r$th row of $A$. Recall that $D_r = \sum_{i=1}^n a_{ir}\phi_i\partial_{\phi_i}=\operatorname{diag}(\phi)A_r$. Therefore, $\operatorname{dim span}\{D_1,D_2,\dots, D_k\} = \operatorname{dim } \operatorname{rowsp}(A) = \operatorname{rank}(A)$, where $\operatorname{rowsp}(A)$ is the row space of $A$. Thus, H\"ormander's condition holds if and only if $D_0\notin \operatorname{span}\{D_1,D_2,\dots, D_k\}$. The condition $D_0\in \operatorname{span}\{D_1,D_2,\dots, D_k\}$ is equivalent to the existence of a vector $c \in \R^k$ such that, for all $i\in [n]$,
    \[\frac{\sum_{j=1}^n \Sigma_{ij}\phi_j}{y(\phi)} - \frac{1}{2}||a_i||^2 = \sum_{r=1}^k c_r a_{ir}.\]
    Substituting $\Sigma_{ij} = a_i\cdot a_j$ gives $a_i\cdot \left(\frac{\sum_{j=1}^n a_j \phi_j}{y(\phi)}\right)-\frac{1}{2}||a_i||^2 = a_i\cdot c$. Rearranging yields $a_i\cdot \left(\frac{\sum_{j=1}^n a_j \phi_j}{y(\phi)}- c\right) = \frac{1}{2}||a_i||^2$. Collecting all coordinates, we can write this system as $A^T z = \frac{1}{2}(||a_1||^2, \dots, ||a_n||^2)^T$. This equation has a solution if and only if the vector $(||a_1||^2, \dots, ||a_n||^2)$ belongs to $\operatorname{colsp}(A^T)$, which is the row space of $A$. Therefore, linear independence of the drift vector field $D_0$ from the diffusion fields requires $(||a_1||^2, \dots, ||a_n||^2) \notin \operatorname{rowsp}(A)$.
\end{proof}
\begin{remark}
    The condition in Theorem \ref{thm_testing} can also be expressed compactly by requiring the following matrix to have full rank:
\begin{equation*}
   \left[
   \begin{array}{ccc|c}
        a_{11} & \cdots & a_{k1} & \|a_1\|^2 \\
        \vdots & \ddots & \vdots & \vdots  \\
        a_{1n} & \cdots & a_{kn} & \|a_n\|^2 
    \end{array}
    \right] .
\end{equation*}
    
\end{remark}
\subsection{The detection case}
When $Q\neq 0$, the generator in the $\phi$-coordinate becomes
\[\calL = \frac{1}{2}\sum_{i,j=1}^n \Sigma_{ij} \phi_i\phi_j \frac{\partial^2}{\partial_{\phi_i}\partial_{\phi_j}} +\frac{1}{Y}\sum_{i,j=1}^n\Sigma_{ij} \phi_i\phi_j \frac{\partial}{\partial_{\phi_j}} + \sum_{j=1}^n\sum_{i=0}^n (q_{ij} - q_{i0}\phi_j)\phi_i\frac{\partial}{\partial_{\phi_j}} \]
where $\Sigma_{ij}:= (\lambda_i - \lambda_0)\cdot (\lambda_j-\lambda_0)$ and $Y:=1+\sum_{i=1}^n \phi_i$. Recall that each $\lambda_i$ is a $k$-dimensional vector. The diffusion vector fields remain unchanged: $D_r = \sum_{j=1}^n a_{jr} \phi_j \partial_{\phi_j}$ for $r \in[k]$, and they still commute. The drift vector field becomes
\[D_0^J = \sum_{j=1}^n \left(\frac{1}{y(\phi)}\sum_{i=1}^n \Sigma_{ij}\phi_i\phi_j - \frac{1}{2}||a_j||^2\phi_j\right)\partial_{\phi_j}+ \sum_{j=1}^n\sum_{i=0}^n (q_{ij} - q_{i0}\phi_j)\phi_i\partial_{\phi_j}.\]
Write $D_0^J = D_0+J$, where $J = \sum_{j=1}^n\sum_{i=0}^n (q_{ij} - q_{i0}\phi_j)\phi_i\partial_{\phi_j}$ represents the contribution from transitions of the hidden Markov chain $\theta$. The derivative of the $j$th component of $J$ with respect to $\phi_m$ is $\partial_{\phi_m} J_j  = q_{mj}- \phi_jq_{m0} - \sum_{i=0}^n\phi_iq_{i0}1_{m=j}$. For any $r\in[k]$,
\[[D_r, D_0^J] = [D_r, D_0]+[D_r,J] = \sum_{s=1}^k \left(\frac{\alpha_{rs}}{y(\phi)} - \frac{\beta_r\beta_s}{y^2(\phi)}\right) D_s+[D_r,J].\]
The $j$th component of $[D_r, J]$ is 
\begin{align*}
    [D_r, J]^j& = \sum_{m=1}^na_{mr}\phi_m\partial_{\phi_m}J_j - \sum_{i=0}^n J_i\partial_{\phi_i} \left(a_{jr}\phi_j\right)\\
    & = \sum_{m=1}^n a_{mr}q_{mj}\phi_m -\sum_{m=1}^n a_{mr}q_{m0}\phi_m\phi_j -a_{jr}\phi_j\sum_{i=0}^n\phi_iq_{i0} - a_{jr}J_j\\
        & =\sum_{m=1}^n (a_{mr}-a_{jr})q_{mj}\phi_m - \phi_j \sum_{m=1}^n a_{mr} q_{m0}\phi_m - a_{jr} q_{0j}.
\end{align*}
Therefore, $[D_r, J]$ is a diagonal vector field whose components are polynomials of degrees $0$, $1$, and $2$ in $\phi$. In particular, it produces directions not generated by the diffusion fields alone.

\begin{lemma}\label{adD}
    For any vector field $U$, let $ad_{D_r}U: = [D_r, U]$. Then $ad_{D_r}$ and $ad_{D_s}$ commute for all $r,s\in[k]$: $ad_{D_r}ad_{D_s} = ad_{D_s}ad_{D_r}$.
\end{lemma}
\begin{proof}
Recall that $[D_r, D_s] = 0$ for all $r,s\in[k]$. By the Jacobi identity,
\[ad_{D_r}ad_{D_s} U - ad_{D_s}ad_{D_r}U = [D_r,[D_s,U]]-[D_s,[D_r,U]] = -[U,[-D_s,-D_r]] = 0.\]
It follows that $ad_{D_r}$ and $ad_{D_s}$ commute.
\end{proof}
Let $\alpha = (\alpha_1,\alpha_2,\dots, \alpha_k)$ be a multiindex with $\alpha_i\in\N$ for each $i\in[k]$. Define the operator $\mathbb{G}^{\alpha}$ by
\[\mathbb{G}^{\alpha} U = ad_{D_1}^{\alpha_1}ad_{D_2}^{\alpha_2}\dots ad_{D_k}^{\alpha_k} U=[D_1,\dots,[D_k,[\dots,[D_1,U],\dots,]]].\]
We now consider diagonal fields $U = \sum_{j=1}^nU_j(\phi)  \partial_{\phi_j}$ whose components $U_j(\phi) =\phi_j (\sum_{m=1}^n b_{mj}\phi_m)+  (\sum_{m=1}^n c_{mj}\phi_m) +d_j$ are polynomials of degree at most $2$. Write $U = U^0+U^1+U^2$, where
\[   U^0 = \sum_{j=1}^n d_j \partial_{\phi_j},\quad
    U^1 = \sum_{j=1}^n (\sum_{m=1}^n c_{mj}\phi_m) \partial_{\phi_j},\quad
    U^2 = \sum_{j=1}^n \phi_j (\sum_{m=1}^n b_{mj}\phi_m) \partial_{\phi_j}.\]
The $j$th component of the bracket $[D_r, U]$ is $ [D_r, U_j]=\sum_{m=1}^na_{mr}\phi_m\partial_{\phi_m}U_j - a_{jr}U_j$. For the constant part, $[D_r, U^0]_j=-d_ja_{jr}$. Applying another diffusion field preserves the form of the constant part: $[D_s,[D_r, U^0]]_j=d_ja_{jr}a_{js}$. Therefore, the constant part of $\mathbb{G}^{\alpha}U$ is $\mathbb{G}^{\alpha}U^0 = (-1)^{\sum_{i=1}^k\alpha_i}\left(\prod_{r=1}^k a_{jr}^{\alpha_r} \right) d_j$. For the linear part,
\begin{align*}
[D_r, U^1]_j&= \sum_{m=1}^na_{mr}\phi_m\partial_{\phi_m}(\sum_{i=1}^nc_{ij}\phi_i) - a_{jr}(\sum_{m=1}^n c_{mj}\phi_m)  \\
& = \sum_{m=1}^n a_{mr}c_{mj}\phi_m -a_{jr}(\sum_{m=1}^n c_{mj}\phi_m)\\
& = \sum_{m=1}^n(a_{mr}-a_{jr} )c_{mj}\phi_m.
\end{align*}
Applying another diffusion field $D_s$ yields
\begin{align*}
 [D_s,[D_r, U^1]]_j& = \sum_{m=1}^n a_{ms}\phi_m (a_{mr}-a_{jr}) c_{mj} -a_{js}\sum_{m=1}^n(a_{mr}-a_{jr} )c_{mj}\phi_m\\
 & =  \sum_{m=1}^n (a_{ms}-a_{js} )(a_{mr}-a_{jr} )c_{mj}\phi_m.
\end{align*}
Therefore, the linear part of $\mathbb{G}^{\alpha}U$ is $\mathbb{G}^{\alpha}U^1 =\sum_{m=1}^n \left(\prod_{r=1}^k (a_{mr}-a_{jr})^{\alpha_r} \right) c_{mj}\phi_m$. Finally, for the quadratic part,
\begin{align*}
[D_r, U^2]_j&= \sum_{m=1}^na_{mr}\phi_m\partial_{\phi_m}(\phi_j (\sum_{i=1}^n b_{ij}\phi_i)) - a_{jr}\phi_j (\sum_{m=1}^n b_{mj}\phi_m) \\
& =\sum_{m=1}^na_{mr}b_{mj} \phi_m\phi_j +a_{jr}\phi_j\sum_{i=1}^n b_{ij}\phi_i - a_{jr}\phi_j (\sum_{m=1}^n b_{mj}\phi_m)\\
& =\phi_j\sum_{m=1}^na_{mr}b_{mj} \phi_m.
\end{align*}
Applying another diffusion field $D_s$ gives $[D_s,[D_r, U^2]]_j=\phi_j \sum_{m=1}^n a_{ms}a_{mr}b_{mj}\phi_m $. Therefore, the quadratic part of $\mathbb{G}^{\alpha}U$ is
$\mathbb{G}^{\alpha}U^2 =\phi_j\sum_{m=1}^n \left(\prod_{r=1}^k a_{mr}^{\alpha_r} \right) b_{mj}\phi_m$. Taking $U=J$, we can write the $j$th component of $\mathbb{G}^{\alpha} J$ as 
\begin{equation}
\label{Galpha}
    (\mathbb{G}^{\alpha} J)_j=  (-1)^{\vert \alpha\vert}\left(\prod_{r=1}^k a_{jr}^{\alpha_r} \right) \left(q_{0j}\right)+\sum_{m=1}^n \left(\prod_{r=1}^k (a_{mr}-a_{jr})^{\alpha_r}  \right) q_{mj}\phi_m-\phi_j\sum_{m=1}^n \left(\prod_{r=1}^k a_{mr}^{\alpha_r} \right) q_{m0}\phi_m ,
\end{equation}
for any $\alpha$ with $\vert \alpha\vert>0$.
\begin{lemma}\label{inclusion}
Let $\mathcal{H}(\phi)$ denote the pointwise span of all iterated Lie brackets of $D_0^J,D_1,\dots,D_k$ that contain at least one diffusion field $D_r$, $r\in[k]$. Then, for every $|\alpha|\ge 1$, we have $\mathbb G^\alpha J(\phi)\in \mathcal{H}(\phi)$ for $\phi \in (0,\infty)^n$. In particular, $\mathbb G^\alpha J(\phi) \in\operatorname{Lie}(D_0^J,D_1,\dots,D_k)(\phi)$.
\end{lemma}
\begin{proof}
Let $\mathcal M_D:=\operatorname{span}_{C^\infty}(D_1,\dots,D_k)$. Since the diffusion vector fields commute, $[D_r,\mathcal M_D]\subseteq \mathcal M_D$. If $M=\sum_{s=1}^k f_s D_s\in\mathcal M_D$, then $[D_r,M]=\sum_{s=1}^k D_r(f_s)D_s\in\mathcal M_D$. By Lemma \ref{testing_bracket}, $[D_r,D_0]\in\mathcal M_D$. Then we have $\mathbb G^{e_r}J=[D_r,J]=[D_r,D_0^J]-[D_r,D_0]$. The first term is an iterated Lie bracket containing $D_r$, while the second belongs to $\mathcal M_D$. Therefore, $\mathbb G^{e_r}J(\phi)\in \mathcal{H}(\phi)$. \par 
Next, we proceed by induction on $|\alpha|$. We show that, for every $|\alpha|\geq 1$, $\mathbb G^\alpha J=X+M$, where $X$ is a finite $C^\infty$-linear combination of iterated Lie brackets containing at least one diffusion field and $M\in \mathcal M_D$. The case $|\alpha|=1$ follows from the preceding argument. Suppose the claim holds for some $\alpha$ with $|\alpha|\ge1$. Then, for any $r\in[k]$,
\[\mathbb G^{\alpha+e_r}J=[D_r,\mathbb G^\alpha J]=[D_r,X]+[D_r,M].\]
Every term in $[D_r,X]$ is again a $C^\infty$-linear combination of iterated brackets containing at least one diffusion field, while $[D_r,M]\in \mathcal M_D$. The induction therefore follows. Hence $\mathbb G^\alpha J(\phi)\in \mathcal{H}(\phi)$ and thus $\mathbb G^\alpha J(\phi)\in\operatorname{Lie}(D_0^J,D_1,\dots,D_k)(\phi)$.
\end{proof}
Recall that the $j$th coordinate of the jump field is
\[J_j(\phi) = q_{0j} +\sum_{m=1}^n q_{mj}\phi_m -\phi_j\sum_{m=1}^n q_{m0 }\phi_m -\phi_j q_{00}.\]
To show that the entire Lie algebra spans $\R^n$, it suffices to identify a subspace that already spans $\R^n$. As shown in Theorem \ref{thm_testing}, when $Q = 0$, the vector fields cannot span more than $k+1$ dimensions. A convenient choice here is the subspace generated by $\mathbb{G}^{\alpha} J$. The next theorem gives a sufficient condition for hypoellipticity based only on these fields, yielding a simple parametric check on the matrix $Q$. The proof constructs $n$ fields that form a basis of $\R^n$ at each point $\phi$. 
\begin{theorem}
\label{thm_suff_1}
    Assume that the drift vectors $a_0, a_1, a_2,\dots, a_n$ are pairwise distinct. Suppose further that, for each coordinate $i\in[n]$, there is an $m\neq i$ such that $q_{mi}>0$, and that $q_{kj}>0$ and $a_k-a_j=a_m-a_i$ imply $j=i$. Then H\"ormander's condition holds for all $\phi\in (0,\infty)^n$.
\end{theorem}
\begin{proof}
    
For each $i\in[n]$, define an operator
\[T^i J: = \left(\prod_{r=1}^k P_r^i(\mathbb{G}^{e_r}) \right)J = \sum_{\alpha\in F^i} C_{\alpha} \mathbb{G}^{\alpha}J \]
for some polynomials $P_1^i,\dots, P_k^i$. In particular, if
$T^iJ=\sum_{\alpha\in F^i} C_\alpha\mathbb G^\alpha J$ and the polynomial defining $T^i$ is chosen so that its constant term vanishes, then $T^iJ
=\sum_{\substack{\alpha\in F^i, |\alpha|\ge1}}C_\alpha\mathbb G^\alpha J$. Then, by Lemma \ref{inclusion}, $T^iJ(\phi)\in \mathcal{H}(\phi) \subset\operatorname{Lie}(D_0^J,D_1,\dots,D_k)(\phi)$. We construct these polynomials and the corresponding $F^i$ and $C_{\alpha}$ in turn. The $j$th coordinate of $T^i J$ is
\[(T^i J)_j(\phi)=  \left(\prod_{r=1}^k P_r^i(-a_{jr}) \right) q_{0j}+\sum_{m=1}^n \left(\prod_{r=1}^k P_r^i(a_{mr}-a_{jr}) \right) q_{mj}\phi_m-\phi_j\sum_{m=1}^n \left(\prod_{r=1}^k P_r^i(a_{mr})\right) q_{m0}\phi_m .
\]
We choose the polynomial values ${P^i_r(\cdot)}$ so that all terms vanish when $j \neq i$, while $(T^iJ)_i$ is nonzero when $j = i$. Consider any $k,j\in\{0,\ldots,n\}$ with $k\neq j$ and $q_{kj}>0$ such that $a_k-a_j\neq a_m-a_i$.
Since these two vectors are different, there exists a coordinate
$r(k,j)\in[k]$ such that $a_{k,r(k,j)}-a_{j,r(k,j)}\neq a_{m,r(k,j)}-a_{i,r(k,j)}$. Assign to the polynomial $P_{r(k,j)}^i$ the factor
\[
\frac{
x-\bigl(a_{k,r(k,j)}-a_{j,r(k,j)}\bigr)
}{
a_{m,r(k,j)}-a_{i,r(k,j)}
-\bigl(a_{k,r(k,j)}-a_{j,r(k,j)}\bigr)
}.
\]
This factor equals one when evaluated at $a_{m,r(k,j)}-a_{i,r(k,j)}$
and zero when evaluated at $a_{k,r(k,j)}-a_{j,r(k,j)}$. Since $m\neq i$ and $a_m\neq a_i$, there is also an $r_0\in[k]$ such that $a_{mr_0}-a_{ir_0}\neq0$. Assign to $P_{r_0}^i$ the additional factor $\frac{x}{a_{mr_0}-a_{ir_0}}$. This factor equals one at the selected state and vanishes at
zero. For each $r\in[k]$, let $P_r^i$ be the product of all factors assigned
to the $r$th coordinate. If no factor is assigned to a given coordinate, set $P_r^i\equiv1$. By construction, $\prod_{r=1}^kP_r^i(a_{mr}-a_{ir})=1$, 
and $\prod_{r=1}^kP_r^i(a_{kr}-a_{jr})=0$ for every $k\to j$ with
$a_k-a_j\neq a_m-a_i$. Moreover, $\prod_{r=1}^kP_r^i(0)=0$. Suppose that $q_{kj}>0$ and the drift differences are the same: $a_k-a_j=a_m-a_i$. By the assumption of the theorem, we must have $j=i$. Hence no term with the selected difference can survive in a coordinate different from $i$. If $j=i$, then $a_k-a_i=a_m-a_i$, and thus $a_k=a_m$. Since $a_0,a_1,\ldots,a_n$ are pairwise distinct, it follows that $k=m$. Thus, the only active transition with the selected drift difference is
the chosen transition $m\to i$. Consequently,
\[
(T^iJ)_j(\phi)=0
\qquad\text{for every }j\neq i.
\]
If $m=0$, then $T^iJ=q_{0i}e_i$, and if $m\in[n]$, then $T^iJ=q_{mi}\phi_m e_i$. In either case, $T^iJ=c_i(\phi)e_i$, where $c_i(\phi)\neq0$ for every $\phi\in(0,\infty)^n$. Repeating this construction for every $i\in[n]$, we obtain $T^1J(\phi),\ldots,T^nJ(\phi)\in \mathcal{H}(\phi)$, where each $T^iJ(\phi)$ is a nonzero multiple of $e_i$. Hence the matrix $\bigl(T^1J(\phi),\ldots,T^nJ(\phi)\bigr)$ is diagonal with nonzero diagonal entries and therefore has rank $n$. Consequently, $\mathcal{H}(\phi)=\mathbb R^n$. Since $\mathcal{H}(\phi)\subseteq
\operatorname{Lie}(D_0^J,D_1,\ldots,D_k)(\phi)$, it follows that $\dim\operatorname{Lie}(D_0^J,D_1,\ldots,D_k)(\phi)=n$, and H\"ormander's condition holds for every $\phi\in(0,\infty)^n$.
\end{proof}
\begin{remark}      
The conditions in Theorem \ref{thm_suff_1} are sufficient but not necessary. Consider a generator matrix $Q$ with exactly one column $i$ such that $q_{0i}=0$ and $q_{mi} = 0$ for all $m\neq i$. If there is a $p$ such that $q_{p0}>0$, we can choose the defining polynomials to take the value $1$ at $a_{pr}$ for every $r$. For every other column $j$, set $P_r^i(a_{mr}) = P_r^i(a_{mr}-a_{jr}) = 0$ for all $m\neq p$. Consequently,
    $T^iJ = q_{p0}\phi_p (\phi_1,\phi_2,\dots, \phi_n)$. The resulting matrix has the form
    \[
    \begin{bmatrix}
    C_1 & &  & q_{p0}\phi_1 & \\
     & \ddots &  & \cdots  &  \\
    & & C_{j-1} & q_{p0}\phi_{j-1}  &  \\
    &  &  & q_{p0}\phi_{j}   & \\
    &  &  & \vdots   &\ddots \\
    &  &  & q_{p0}\phi_{j}   & C_n 
\end{bmatrix},\]
and its determinant is nonzero. Therefore, H\"ormander's condition still holds.\par
If there is more than one such column, however, the resulting vectors are collinear, and this construction fails. This observation motivates the alternative sufficient condition in Proposition \ref{thm_suff_2}.\par
Failure of this construction does not imply failure of H\"ormander's condition, because other constructions may work. Importantly, the theorem provides an easily checked but deliberately conservative parametric condition: it considers only the subalgebra generated by iterating the field $J$ with the diffusion and drift fields. Other constructions may also exploit the original directions in $\phi$-space.
\end{remark}
The classical quickest disorder detection problem is a special case covered by Theorem \ref{thm_suff_1}. It corresponds to the setting in which, for each coordinate $i\in[n]$, $q_{0i}>0$ and $q_{mj} = 0$ for all $m\neq j$. In other words, once the drift changes, it remains constant. We refer readers to Example \ref{eg2} for a formulation in this setting.
Theorem \ref{thm_suff_1} gives a simple parametric sufficient condition on the infinitesimal generator $Q$ under which the Lie algebra spans $\R^n$. The following proposition gives another sufficient condition for hypoellipticity when an assumption of Theorem \ref{thm_suff_1} fails; that is, when there is a $j\in [n]$ with no positive incoming intensity, so $q_{mj} = 0$ for all $m\neq j$. In this case, the Lie algebra has limited ability to expand beyond $\R^k$, but it may still span $\R^n$ under certain conditions.
\begin{proposition}
    \label{thm_suff_2}
    Assume that, for every $j\in [n]$, $q_{mj} = 0$ for all $m\in\bar n$ with $m\neq j$, and that there is at least one $j\in[n]$ such that $q_{j0}>0$. Define the augmented matrix
    \begin{equation*}
   \tilde A: = \left[
   \begin{array}{ccc|c|c}
        a_{11} & \cdots & a_{k1} & q_{10}+\frac{1}{2}\|a_1\|^2 & 1 \\
        \vdots & \ddots & \vdots & \vdots & \vdots \\
        a_{1n} & \cdots & a_{kn} & q_{n0}+\frac{1}{2}\|a_n\|^2 & 1 
    \end{array}
    \right] \in \R^{n\times (k+2)}.
\end{equation*}
At any interior point $\phi\in (0,\infty)^n$, we have
    $\operatorname {dim Lie}\{ D_0^J, D_1,D_2,\dots, D_k\} =min\{\operatorname {rank}(\tilde A), n\}$.
    \end{proposition}
    \begin{proof}
        Recall that, for any multiindex $\alpha$,
        \[\mathbb{G}^{\alpha} J = \phi\sum_{m=1}^n \left(\prod_{r=1}^k P(a_{mr})\right) q_{m0}\phi_m \]
        for some polynomial $P$. For a chosen $m$ with $q_{m0}>0$, we can choose $P$ such that $P(a_{mr})=1$ for all $r$ and, for any other $i\neq m$, $P(a_{ir})=1$ for some $r$. This construction gives the vector $q_{m0}\phi_m \phi$, which can be normalized to $\phi$. As in Theorem \ref{thm_testing}, the matrix of basis coefficients is $\operatorname {diag}(\phi_1,\dots, \phi_n) \tilde A$.
    Since $\operatorname {diag}(\phi_1,\dots, \phi_n)$ is invertible, $\operatorname {diag}(\phi_1,\dots, \phi_n) \tilde A$ and $\tilde A$ have the same rank. In particular, if $\operatorname {rank}(\tilde A) = n$, H\"ormander's condition holds at every interior point $\phi$. It fails automatically, however, when $n>k+2$.
    \end{proof}
The three blocks in $\tilde A$ arise from three different vector fields. The diffusion fields generate $\operatorname {diag}(\phi) \operatorname {colsp}(A^T)$. All positive-order fields $\mathbb{G}^{\alpha}J$ are proportional to $\operatorname {diag}(\phi) \mathrm 1$ and thus provide the column $1$. Finally, the field $D_0^J$ contributes the column $(q_{i0}+\frac{1}{2}\|a_i\|^2)_{i\in[n]}$. This decomposition also explains how H\"ormander's condition differs from its parabolic counterpart, as discussed in Section \ref{sec_consequence}.\par
We refer readers to Example \ref{eg3} for a formulation in this setting. The sufficient conditions in Theorem \ref{thm_suff_1} and Proposition \ref{thm_suff_2} impose graph conditions that differ from irreducibility. The matrix $Q$ can be reducible, so some states may be unreachable from others, while H\"ormander's condition still holds and propagates noise in every direction.
\section{Parabolic and joint H\"ormander conditions and consequences of hypoellipticity}
\label{sec_consequence}
The previous section concerns the stationary operator $\calL_{\phi}$ of the posterior likelihood process. Depending on the underlying stopping or control problem, three closely related operators also arise. If the value function depends on time, the relevant operator is the parabolic posterior operator $-\partial_t+\calL_{\phi}$. If the payoff depends directly on the observation process $X$, the state variable becomes $(\phi,x)$, and the relevant operator is the stationary joint operator $\calL_{\phi,x}$. Finally, if the problem depends on both time and the observation process, the corresponding operator is $-\partial_t+\calL_{\phi,x}$. The H\"ormander conditions for these four operators are not equivalent. The purpose of this section is to identify precisely how the four conditions are related in general, in the testing case, and under the two structural assumptions of Theorem \ref{thm_suff_1} and Proposition \ref{thm_suff_2}. In particular, the assumptions of Theorem \ref{thm_suff_1} are sufficient for all four operators to be hypoelliptic. \par
To analyze the parabolic H\"ormander condition, recall that $\mathcal H(\phi)$ denotes the pointwise span of all iterated Lie brackets containing at least one diffusion field. Then $\operatorname{Lie}\{D_0^J,D_1,\dots,D_k\}(\phi)=\operatorname{span}\{D_0^J(\phi)\}+\mathcal H(\phi)$. The distinction between the stationary and parabolic H\"ormander conditions can then be stated explicitly.
\begin{proposition}[Stationary and parabolic H\"ormander conditions]
Define the backward time-space operator $\bar D_0^J : = -\partial_t+D_0^J$. Then $\operatorname{dim Lie}\{\bar D_0^J, D_1,D_2,\dots, D_k\}(t,\phi) =1 + \operatorname{dim} \ \mathcal{H}(\phi)$. Thus $\mathcal L_{\phi}$ satisfies H\"ormander's condition at $\phi$ if and only if $\operatorname{span}\{D_0^J(\phi)\} + \mathcal{H}(\phi) = \mathbb R^n$, whereas $-\partial_t+\mathcal L_{\phi}$ satisfies H\"ormander's condition at $(t,\phi)$ if and only if $\mathcal{H}(\phi)=\mathbb R^n$. Hence the parabolic H\"ormander condition implies the stationary H\"ormander condition, but the converse need not hold. If the parabolic H\"ormander condition holds, the process $(\Phi_t)_{t\geq 0}$ is a strong Feller process on $(0,\infty)^n$.
\end{proposition}
\begin{proof}
     Since all the spatial vector fields are independent of $t$, $[\partial_t,D_0^J]=0$ and $[\partial_t,D_r]=0$ for $r\in[k]$. Therefore, $[\bar D_0^J,D_r] = [D_0^J,D_r]$. More generally, every iterated bracket that involves $\bar D_0^J,D_1,\dots,D_k$ and contains at least one diffusion field has zero time component, and its spatial component is exactly the bracket obtained by replacing $\bar D_0^J$ with $D_0^J$. Thus, all such brackets span precisely $ \{0\}\times  \mathcal{H}(\phi)$. On the other hand, since $[\bar D_0^J,\bar D_0^J]=0$, it follows that
     \[\dim \operatorname{Lie} \{\bar D_0^J,D_1,\dots,D_k\}(t,\phi) = 1+\dim  \mathcal{H}(\phi).\]
     The stated characterizations of the two H\"ormander conditions follow immediately. The strong Feller property of $\Phi$ on $(0,\infty)^n$ follows from an argument similar to that used in Proposition 4 of \cite{peskir_detection}.
\end{proof}
The strong Feller property here is an interior statement because the H\"ormander conditions above are verified only on $(0,\infty)^n$. Whether this property extends to processes started on the boundary depends on the transition structure of $Q$ and on whether hypoellipticity can be propagated to those points. For example, if $Q$ is irreducible, the process $\Pi$ is strong Feller on the entire closed simplex $P_{n+1}$.
\begin{corollary} \label{parabolic_cor} The following statements describe the parabolic H\"ormander condition:
\begin{enumerate}
    \item \textbf{Testing case.} Assume $Q=0$. Then $\dim \mathcal{H}(\phi)=\operatorname{rank}(A)\leq k<n$. Consequently, the parabolic posterior operator does not satisfy H\"ormander's condition.
    \item \textbf{Under the assumptions of Theorem \ref{thm_suff_1}.} If the assumptions of Theorem \ref{thm_suff_1} hold, then
    \[\operatorname{dim Lie}\{\bar D_0^J, D_1,D_2,\dots, D_k\} =n+1,\]
    and thus the parabolic H\"ormander condition holds.
    \item \textbf{Under the assumptions of Proposition \ref{thm_suff_2}.}
If the assumptions of Proposition \ref{thm_suff_2} hold, then $\mathcal{H}(\phi) = \operatorname{diag}(\phi)
\operatorname{colsp}[A^T\mid\mathbf 1]$. In particular, when $n=k+1$, the parabolic H\"ormander condition holds if $\operatorname{rank}[A^T\mid\mathbf 1]=n$; it fails automatically if $n>k+1$.
\end{enumerate}
\end{corollary}
\begin{proof}
In the testing case, Lemma \ref{testing_bracket} implies that every Lie bracket containing at least one diffusion field remains in the span of $D_1,\dots,D_k$. Hence $\mathcal{H} (\phi)= \operatorname{span}\{D_1(\phi),\dots,D_k(\phi)\} = \operatorname{diag}(\phi)\operatorname{colsp}(A^T)$, so $\dim \mathcal{H}(\phi)\leq k<n$. The parabolic condition therefore fails by the preceding proposition. For the second case, the proof of Theorem \ref{thm_suff_1} constructs, for each $i\in[n]$, a field $T^iJ=c_i(\phi)e_i$ with $c_i(\phi)\neq0$, where $T^iJ$ is a linear combination of fields $\mathbb G^\alpha J$ with $|\alpha|\geq1$. Hence $T^iJ\in \mathcal{H}(\phi)$ for every $i$. The conclusion again follows from the preceding proposition. Example \ref{eg2} is covered by this case. For the third case, the diffusion fields generate $ \operatorname{diag}(\phi)\operatorname{colsp}(A^T)$, while the positive-order fields $\mathbb G^\alpha J$ generate the additional direction $\operatorname{diag}(\phi)\mathbf 1$. Thus $\mathcal{H}(\phi) = \operatorname{diag}(\phi)
\operatorname{colsp}[A^T\mid\mathbf 1]$, and the parabolic characterization follows from the previous proposition.
\end{proof}
In some stopping problems with partial information, the payoff functions depend on both the unknown state $\theta$ and the observation process $X$. Consequently, the value function depends on $(\pi,x)$. The following proposition relates hypoellipticity in the $\phi$-coordinates to hypoellipticity in the joint $(\pi,x)$-coordinates. Recall that the $k$-dimensional observation process satisfies
\[dX_t =  \frac{1}{Y_t}\left(\lambda_0+ \sum_{i=1}^n\lambda_i \Phi_t^i \right) dt +d \tilde W_t =: \tilde \lambda(\Phi_t)dt+ d\tilde W_t. \]
For $r\in[k]$, the $r$th component of $\tilde \lambda(\Phi_t)$ is $\tilde \lambda_r(\Phi_t) = \frac{1}{Y_t}\left(\lambda_{0r}+ \sum_{i=1}^n\lambda_{ir} \Phi_t^i \right)$.
\begin{align*}
  \calL_{\phi,x} = &\frac{1}{2}\sum_{i,j=1}^n \Sigma_{ij} \phi_i\phi_j \frac{\partial^2}{\partial_{\phi_i}\partial_{\phi_j}} +\frac{1}{Y}\sum_{i,j=1}^n\Sigma_{ij} \phi_i\phi_j \frac{\partial}{\partial_{\phi_j}} + \sum_{j=1}^n\sum_{i=0}^n (q_{ij} - q_{i0}\phi_j)\phi_i\frac{\partial}{\partial_{\phi_j}}\\
  & + \frac{1}{2} \sum_{j=1}^k \frac{\partial^2}{\partial_{x_j^2}}+ \sum_{j=1}^k \tilde \lambda_j(\phi)\frac{\partial}{\partial_{x_j}}+\sum_{i=1}^n\sum_{j=1}^ka_{ij} \phi_i\frac{\partial^2}{\partial_{\phi_i}\partial_{x_j}}.
\end{align*}
\begin{proposition}
\label{prop_x}
If $\calL_{\phi,x}$ satisfies H\"ormander's condition, then $\calL_{\phi}$ also satisfies it. Conversely, suppose that $\calL_{\phi}$ satisfies H\"ormander's condition and that there exist operators $\tilde T^1, \dots , \tilde T^n$ such that $\tilde T^i J = \sum_{|\alpha|\geq 1} c_{\alpha}^i \mathbb{G}^{\alpha} J$ for $i\in [n]$, with $\tilde T^1 J(\phi), \dots , \tilde T^n J(\phi)$ linearly independent. Then $\calL_{\phi,x}$ satisfies H\"ormander's condition.
\end{proposition}

\begin{proof}
For each $r\in[k]$, define the vector field $\tilde D_r : = D_r + \partial_{x_r}$, and let $\hat D_0 = D_0^J+D^X = D_0^J +\tilde \lambda(\phi)\cdot \nabla_x$. Observe that $\calL_{\phi,x} = \hat D_0 + \frac{1}{2}\sum_{r=1}^k \tilde D_r^2$. Since the coefficients of $D_r$ depend only on $\phi$, $[\tilde D_r, \tilde D_s] = [D_r,D_s] = 0$. Moreover, the coefficients of $D^X$ depend only on $\phi$, so $[\hat D_0, \tilde D_r] = [D_0^J, D_r] - \sum_{s=1}^k D_r(\tilde \lambda_s)\partial_{x_s}$. We then have $[\tilde D_u,[\hat D_0,\tilde D_r]] = [ D_u,[D_0^J,D_r]] - \sum_{s=1}^k D_u D_r(\tilde \lambda_s)\partial_{x_s}$. Similarly, every iterated bracket of $\hat D_0,\tilde D_1,\dots, \tilde D_k$ has the form
\[[D_0^J,[D_u,\dots, [D_0^J,D_r]\dots]]+\sum_{s=1}^k F_s(\phi)\partial_{x_s}.\]
Consider the projection onto the $\phi$-coordinates that sends $\hat D_0$ to $D_0^J$ and $\tilde D_r$ to $D_r$. This projection maps each iterated bracket of $\hat D_0,\tilde D_1,\dots, \tilde D_k$ to the corresponding bracket of $D_0^J,D_1,\dots , D_k$. It follows that if the joint Lie algebra has dimension $n+k$, then the posterior-likelihood Lie algebra has dimension $n$.\par
Conversely, suppose that the additional assumption holds. Then
\[[\tilde D_r,\hat D_0] = [D_r, D_0^J]+\sum_{s=1}^k D_r(\tilde \lambda_s)\partial_{x_s} = [D_r, J] + \sum_{s=1}^k D_r(\tilde \lambda_s)\tilde D_s.\]
Thus, at every $(\phi,x)$, $[\tilde D_r,\hat D_0] - \sum_{s=1}^k D_r(\tilde \lambda_s)\tilde D_s = \mathbb{G}^{e_r}J\in\Span \operatorname{Lie }\{\hat D_0,\tilde D_1,\dots, \tilde D_k\}(\phi,x)$. Since $\mathbb{G}^{\alpha}J$ depends only on $\phi$, $[\tilde D_s, \mathbb{G}^{\alpha}J] = \mathbb{G}^{\alpha+e_s}J$. It follows that $\mathbb{G}^{\alpha}J(\phi) \in  \Span\operatorname{Lie }\{\hat D_0,\tilde D_1,\dots, \tilde D_k\}(\phi,x)$ for every $|\alpha|\geq 1$. By the additional assumption, $\tilde T^i J(\phi) \in  \Span\operatorname{Lie }\{\hat D_0,\tilde D_1,\dots, \tilde D_k\}(\phi,x)$ for every $i\in [n]$, and these fields span all $n$ directions in the $\phi$-coordinates. Therefore, for every $r\in [k]$, there are coefficients $b_{ri}(\phi)$ such that $D_r(\phi) = \sum_{i=1}^n b_{ri}(\phi)\tilde T^iJ(\phi)$ and $\tilde D_r(\phi,x) - \sum_{i=1}^n b_{ri}(\phi)\tilde T^iJ(\phi) = \partial_{x_r}$. The joint Lie algebra therefore contains $n$ directions in $\phi$ and $k$ directions in $x$, so $\operatorname{dim Lie }\{\hat D_0, \tilde D_1,\dots, \tilde D_k\}(\phi,x) = n+k$.
\end{proof}
The proposition therefore covers many cases in which the running payoff $h$, the immediate payoff $g$, or both depend on $(\pi,x)$.
\begin{corollary}\label{joint_cor}
The following statements describe the joint H\"ormander condition:
\begin{enumerate}
    \item \textbf{Testing case.} Assume $Q=0$. Then $\mathcal L_{\phi,x}$ does not satisfy H\"ormander's condition in the degenerate case $k<n$.
    \item \textbf{Under the assumptions of Theorem \ref{thm_suff_1}.} If the assumptions of Theorem \ref{thm_suff_1} hold, then $\calL_{\phi,x}$ satisfies H\"ormander's condition.
    \item \textbf{Under the assumptions of Proposition \ref{thm_suff_2}.}
Assume that the hypotheses of Proposition \ref{thm_suff_2} hold. If $n\geq3$, then $\mathcal L_{\phi,x}$ does not satisfy H\"ormander's condition. If $n=2$ (and hence $k=1$), then $\mathcal L_{\phi,x}$ satisfies H\"ormander's condition if and only if $q_{10}+\frac{1}{2}||a_1||^2+\lambda_0\cdot a_1\neq q_{20}+\frac{1}{2}||a_2||^2+\lambda_0\cdot a_2$.
\end{enumerate}
\end{corollary}
\begin{proof}
In the testing case, working in logarithmic coordinates, we have $\tilde D_r (\log \phi_i - a_i\cdot x) =0$ and $\hat D_0(\log \phi_i - a_i\cdot x) = -\frac12\|a_i\|^2-\lambda_0\cdot a_i$. Thus, the Lie algebra generated by $\hat D_0,\tilde D_1,\dots,\tilde D_k$ has dimension at most $k+1$. Since the joint state space has dimension $n+k>k+1$, H\"ormander's condition fails. The second claim follows because the construction of $\tilde T^i J$ directly satisfies the additional condition in Proposition \ref{prop_x}. Under the assumptions of Proposition \ref{thm_suff_2}, however, the situation is different. In logarithmic coordinates, $\tilde D_r(\log \phi_i - a_i\cdot x) = 0$ and $\hat D_0 (\log \phi_i - a_i\cdot x) = f(\phi ) -(q_{i0}+\frac{1}{2}||a_i||^2+\lambda_0\cdot a_i)$. For $n\geq 3$, we can choose $c\neq 0$ such that $\sum_{i=1}^n c_i = 0$ and $\sum_{i=1}^n c_i\left(q_{i0}+\frac{1}{2}\|a_i\|^2+\lambda_0\cdot a_i\right)=0$. Consequently, $\tilde D_r
\left(\sum_{i=1}^n c_i(\log\phi_i-a_i\cdot x)
\right)=0$ and $\hat D_0\left(\sum_{i=1}^n c_i(\log\phi_i-a_i\cdot x)\right)=0$, so the joint H\"ormander condition fails. When $n=2$, it holds if and only if $q_{10}+\frac{1}{2}||a_1||^2+\lambda_0\cdot a_1\neq q_{20}+\frac{1}{2}||a_2||^2+\lambda_0\cdot a_2$.
\end{proof}
Finally, we consider the parabolic joint operator $-\partial_t+\calL_{\phi,x}$. The proof of the following result is parallel to that of Corollary \ref{parabolic_cor}.
\begin{proposition}
    \label{prop_parabolic_joint}
    Consider the parabolic joint operator $-\partial_t+\calL_{\phi,x} = -\partial_t+\hat D_0 + \frac{1}{2}\sum_{r=1}^k \tilde D_r^2$. Then
    \begin{enumerate}
        \item In the testing case or under the assumptions of Proposition \ref{thm_suff_2}, H\"ormander's condition fails.
        \item Under the assumptions of Theorem \ref{thm_suff_1}, H\"ormander's condition holds.
    \end{enumerate}
\end{proposition}
\begin{proof}
All the assertions are pointwise on the interior time-space domain $\R\times(0,\infty)^n\times\R^k$. In the testing case, set $I_i(\phi,x):=\log\phi_i-a_i\cdot x$ and $b_i:=\frac12\|a_i\|^2+\lambda_0\cdot a_i$. We have $\tilde D_r I_i=0$ for every $r\in[k]$ and $\hat D_0 I_i=-b_i$. Hence
\[
F_i(t,\phi,x):=I_i(\phi,x)-b_it
\]
is a nonconstant common first integral of $-\partial_t+\hat D_0,\tilde D_1,\dots,\tilde D_k$. Indeed, $(-\partial_t+\hat D_0)F_i=b_i-b_i=0$, and $\tilde D_rF_i=0$ for every $r$. Every iterated Lie bracket therefore annihilates $F_i$, so the Lie algebra is tangent to the proper level sets of $F_i$ and cannot have full rank. Thus, the parabolic H\"ormander condition fails in the testing case.\par
Under Proposition \ref{thm_suff_2}, let $b_i:=q_{i0}+\frac12\|a_i\|^2+\lambda_0\cdot a_i$. As shown above, $\tilde D_r I_i=0$ and $\hat D_0 I_i=f(\phi)-b_i$, where $f$ does not depend on $i$. Choose $c\neq0$ with $c\cdot\mathbf 1=0$ and define
\[
F_c(t,\phi,x):=\sum_{i=1}^n c_iI_i(\phi,x)-t\sum_{i=1}^n c_ib_i.
\]
Then $\tilde D_rF_c=0$ for every $r$, while
\[
(-\partial_t+\hat D_0)F_c
=\sum_{i=1}^n c_ib_i+f(\phi)\sum_{i=1}^n c_i-\sum_{i=1}^n c_ib_i=0.
\]
Since $c\neq0$, $F_c$ is nonconstant. The same first-integral argument shows that the parabolic H\"ormander condition fails.\par
Finally, under Theorem \ref{thm_suff_1}, the positive-order brackets constructed above span all posterior directions, and Proposition \ref{prop_x} then yields all $n+k$ spatial directions. These directions arise from brackets containing at least one diffusion field and therefore have zero time component. Since $-\partial_t+\hat D_0$ has time component $-1$, it supplies one additional independent direction. Hence the parabolic joint Lie algebra has rank $n+k+1$, and H\"ormander's condition holds.
\end{proof}
We consider the following variational inequality associated with the stopping problem \eqref{value_fcn} in $\mathring{P}_{n+1}$: 
\begin{equation}
\label{VI}
    \min\{r u-\calL_{\pi} u -h,u-g \} = 0.
\end{equation}
As argued above, a posterior probability process starting in $\mathring{P}_{n+1}$ remains there. Therefore, the local regularity analysis below does not involve the boundary $\partial P_{n+1}$. We do not address a global viscosity characterization or the corresponding boundary conditions. For the remainder of this section, we assume that $V$ is continuous in $\mathring P_{n+1}$.
Define the continuation region
\[\C: = \{\pi\in \mathring{P}_{n+1}: V(\pi)>g(\pi)\}\]
and the stopping region
\[\D: = \{\pi\in\mathring{P}_{n+1}: V(\pi)=g(\pi)\}.\]
When an optimal first-entry rule is available under appropriate standard optimal-stopping assumptions, $\D$ is the corresponding stopping region. By the continuity of $V-g$, the continuation region is open and the stopping region is closed. Let $U\Subset\C$ and let
$\tau_U:=\inf\{t\geq0:\Pi_t\notin U\}$. The localized dynamic programming principle implies that, for $\pi\in U$, $V$ satisfies
\begin{equation}
    \label{DPP}
V(\pi) = \E_\pi\left[\int_0^{\tau_U}
e^{-\int_0^t r(\Pi_s)ds}h(\Pi_t)dt+
e^{-\int_0^{\tau_U}r(\Pi_s)ds}V(\Pi_{\tau_U})\right].\end{equation}
Consequently, $V$ is a weak solution of $(r-\calL_\pi)V=h$ in $\C$. The continuation-region analysis above is carried out entirely in $\mathring P_{n+1}$. Boundary points $\pi\in\partial P_{n+1}$ must be treated separately. As observed in Section \ref{sec_change}, a process starting in the interior remains there for all $t\geq 0$ with probability $1$. Thus, the boundary $\partial {P}_{n+1}$ is unattainable from the interior. For a process starting on the boundary, fix an index $j\in\bar{n}$ and consider the face $F_j:=\{\pi\in {P}_{n+1}, \enskip \pi_j = 0\}$. When $\Pi_0 = \pi \in F_j$, the diffusion coefficient in the $j$th direction vanishes because $\pi_j(\lambda_j - \bar\lambda_0) = 0$, while the drift coefficient satisfies $\sum_{i\neq j} q_{ij}\pi_i\geq 0$. In the testing case ($Q=0$), $\sum_{i\neq j} q_{ij}\pi_i=0$, so $\Pi$ does not move in the $j$th direction and $\Pi_t^j\in F_j$ for all $t\geq 0$. Consequently, the boundary of the simplex serves as a natural boundary for the variational inequality \eqref{VI}, and the optimal stopping problem reduces to an $(n-1)$-dimensional problem on $F_j$. The value function \eqref{value_fcn} can instead be characterized on the lower-dimensional domain $F_j$.\par
In the detection case ($Q\neq 0$), if $\sum_{i\neq j} q_{ij}\pi_i>0$ for $\pi \in F_j$, the process acquires an immediate positive drift in the $j$th coordinate and leaves the face $F_j$. Thus, $F_j$ is an entrance boundary in the $j$th direction. Although the stopping problem \eqref{value_fcn} defines $V$ for initial states on the closed simplex, identifying its boundary values with limits of the interior value function requires a separate boundary-continuity argument. We do not assert such continuity here because it is not needed for the interior variational inequality or the interior regularity result below.
At boundary points where $\sum_{i\neq j} q_{ij}\pi_i=0$, the process initially diffuses on the lower-dimensional face $F_j$ until it reaches an entrance boundary point. In the special case $q_{ij} = 0$ for all $i$, the process remains in $F_j$, and the problem reduces to a lower-dimensional one, as in the testing case.\par 
Under suitable regularity assumptions, hypoellipticity of $\calL_{\pi}$ implies that the value function $V$ is smooth in the continuation region.
\begin{theorem}
\label{smooth_in_C}
Consider the value function defined in \eqref{value_fcn}. Assume that $r,h\in C^{\infty}$ on $\C$. If H\"ormander's condition holds on $\C$, then the value function is $C^{\infty}$ on $\C$.
\end{theorem} 
The vector fields $D_0^J, D_1,D_2,\dots, D_k$ have $C^{\infty}$ coefficients, so the conclusion follows directly from Corollary $7$ of \cite{peskir_weak}. Consequently, the solution is classical in the continuation region, where It\^o's formula may be applied. \par 
There is a substantial literature on the regularity of equations and obstacle problems associated with H\"ormander-type operators; see, for example, \cite{MF_regularity, LL_schauder, FGN}. In particular, Schauder-type estimates
can be obtained in H\"older spaces adapted to the underlying H\"ormander geometry under suitable additional assumptions. This suggests that some interior regularity may hold under conditions weaker than those in Theorem \ref{smooth_in_C}. We do not discuss such estimates further because they require additional geometric and regularity assumptions.\par
For stopping problems, $C^{\infty}$ regularity holds only in the open continuation set $\C$; it does not extend to the stopping boundary or the boundary of the full domain. In particular, our results do not establish global $C^1$ regularity across $\partial \C$ and therefore do not imply smooth fit. Existing results address global $C^1$ regularity of value functions for optimal stopping problems. In our setting, suppose that the underlying process $\Pi$ is strong Feller, the assumptions of Theorem \ref{smooth_in_C} hold, $g$ is $C^1$, and every point of $\partial \C$ is probabilistically regular. Then the value function is $C^1$ everywhere. For details, see Theorem $8$ in \cite{DAP}.
\section{Examples}\label{sec_example}
\begin{example}(Testing three possible drifts of a one-dimensional Brownian motion)\\
\label{eg1}
    Consider the testing case with $k=1$ and $n=2$. We observe a one-dimensional Brownian motion whose constant drift can take one of three values: $\lambda_0=0$, $\lambda_1 = 1$, or $\lambda_2=2$. We are interested in testing the following three hypotheses about the drift:
\[H_0: \theta = 0, \enskip H_1: \theta = 1,  \enskip H_2: \theta = 2.\]
The observer pays a constant cost $c>0$ per unit of time and incurs a penalty for an incorrect decision. The observer therefore balances the benefit of greater accuracy against the cost of additional observation. Let the decision rule $d\in\{0,1,2\}$ represent the stopper's estimate of the true value of $\theta$. The objective is to minimize the observation cost while making an accurate decision:
\[V = \inf_{\tau,d}\E[c\tau + \sum_{i = 0}^2 b_{i}1_{\{\theta = i, d\neq i\}}],\]
where $b_i>0$ denotes the penalty for deciding $d\neq i$ when $\theta = i$. The problem can equivalently be written as
\[V(\pi) = \inf_{\tau,d}\E[c\tau + \min\{b_1\Pi_{\tau}^1+b_2\Pi_{\tau}^2,b_0\Pi_{\tau}^0+b_2\Pi_{\tau}^2,b_0\Pi_{\tau}^0+b_1\Pi_{\tau}^1\}]. \]
In this setting, $A = (1,2)$ and $\operatorname{rank}(A)  = 1$. Moreover, the vector $(\|a_1\|^2,\|a_2\|^2 ) = (1,4)$ does not belong to the row space of $A$. Thus, H\"ormander's condition is satisfied, and Theorem \ref{thm_testing} implies hypoellipticity. Since the running payoff is smooth, the value function is smooth in the continuation region by Theorem \ref{smooth_in_C}. \par
In this problem, the posterior likelihood process $\Phi$ can be written as a function of $t$ and the observation process $X$. Define
$F(t, X_t): =  \frac{\lambda_1\phi_t^1+\lambda_2\phi_t^2}{1+\phi_t^1+\phi_t^2}$ and write $dX_t = F(t,X_t)dt+ d\tilde{W}_t$. One can then formulate an equivalent parabolic problem with 
\[\calL = \partial_t+F(t,x)\partial_x+\frac{1}{2}\partial_{xx},\]
and the payoff function depends on $(t,x)$ and takes a different form in each of the three regions. Standard parabolic tools can then be applied to establish smoothness of the value function. We can also verify hypoellipticity directly by setting $D_0 = \partial_t+F(t,x)\partial_x$ and $D_1 = \partial_x$.
Then $[D_0,D_1]=  -\partial_x (F(t,x)) \partial_x$, so $\dim\left(\operatorname{Lie}\left(D_0, D_1\right)\right) = 2$. Therefore, hypoellipticity is immediate in this one-dimensional example. This example appears in Section $11$ of \cite{caffarelli}; it was also studied as a parabolic problem in the $(t,x)$-coordinates in \cite{zhitlukhin}. However, both the parabolic H\"ormander condition and the joint H\"ormander condition for $\calL_{\phi,x}$ fail in this example. Note that the parabolic operator studied in \cite{zhitlukhin} differs from $-\partial_t+\calL_{\phi}$: it is a lower-dimensional reparametrization that exploits the deterministic relationship between $(t,x)$ and the posterior. Thus, there is no contradiction with Section \ref{sec_consequence}.
\end{example}
As a detection counterpart to Example \ref{eg1}, consider $n=2$, $k=1$, and $\lambda_0=0$, $\lambda_1 = 1$, $\lambda_2=2$. Since $\rank[A^T\vert 1] = 2$, both the stationary and parabolic H\"ormander conditions hold. The joint H\"ormander condition, however, may hold or fail depending on the relationship between the transition rates and the observation drifts. For example, the joint condition holds if $q_{10} = q_{20}=1$ but fails if $q_{10}=\frac{5}{2}$ and $q_{20} = 1$. 
\begin{example}(Detecting the change time in multiple coordinates)\\
\label{eg2}
     Suppose that we observe a standard $N$-dimensional Brownian motion $X = (X^1,\dots, X^N)$ whose drifts are initially zero. Assume that there is a random time $\eta$ such that
    \[\P(\eta=0)=\pi,\quad \P(\eta>t\vert\eta>0)=e^{-\lambda t}, \quad t\geq 0.\]
    On the event $\{\eta>0\}$, exactly $K$ of the $N$ processes permanently acquire a drift $\mu\neq 0$ at time $\eta$. The objective is to detect this change time as accurately as possible. This yields a detection problem with $k=N$ observation processes, $n=\binom{N}{K}$ possible post-change drift vectors, and one zero-drift vector representing no change. \par
    Let $S_1, \dots, S_{\binom{N}{K}}\subset \{1,2,\dots, N\}$ be the $\binom{N}{K}$ distinct subsets of size $K$, representing the coordinates with positive drifts in each scenario. The matrix $A = (a_1, a_2,\dots, a_{\binom{N}{K}})$ then has pairwise distinct columns $a_j\in \mathbb{R}^N$, where
    \begin{align}
        (a_j)_r = \begin{cases} \mu, \quad r\in S_j,\\0, \quad r\notin S_j.
    \end{cases}
    \end{align}
    The generator matrix $Q$, of size $(\binom{N}{K}+1)\times (\binom{N}{K}+1)$, has nonzero entries only in its first row because the drifts change at most once:
    \begin{align}
    Q = \begin{bmatrix} -\lambda & \lambda/\binom{N}{K}  &\lambda/\binom{N}{K} &\dots& \lambda/\binom{N}{K} \\
    0&0&0&\dots& 0\\
   \vdots & \vdots & \vdots & \ddots & \vdots \\
0&0&0&\dots& 0\
    \end{bmatrix}.\end{align}
   Since all drifts are distinct and all $q_{0j}$ are strictly positive, Theorem \ref{thm_suff_1} implies H\"ormander's condition. This problem is studied in detail in \cite{peskir_detection}. The coordinates may also acquire different drifts $\mu_i'$ with different transition rates $\lambda_i'$ without affecting hypoellipticity. Moreover, because the example satisfies the assumptions of Theorem \ref{thm_suff_1}, the parabolic H\"ormander condition and the joint H\"ormander condition for $\calL_{\phi,x}$ hold automatically.
\end{example}
\begin{example}(Sequential tracking with regime switching)\\
\label{eg3}
Suppose that we observe a stochastic process $X$ driven by a $k$-dimensional Brownian motion and having a drift determined by an unobservable Markov process $\theta$ with $n+1$ states:
\[X_t = \int_{0}^t \lambda_{\theta_s}ds+W_t.\] 
Assume that the infinitesimal generator has the form
    \begin{align}
    Q = \begin{bmatrix} -\sum_{i=1}^n q_i & q_1 &q_2&\dots& q_n\\
    p_1&-p_1&0&\dots& 0\\
   \vdots & \vdots & \vdots & \ddots & \vdots \\
p_n&0&0&\dots& -p_n\
    \end{bmatrix}\end{align}
    with $p_i, q_j>0$ for all $i, j\in [n]$. We consider the problem of minimizing the total cost
    \[\inf_{d}\E[\int_0^T e^{-rt} c(d_t,\theta_t)dt],\]
    where $d_t$ takes values in $[n]$ and represents the current estimate of the hidden regime, while the running cost $c(d_t,\theta_t)$ penalizes both delays in detecting switches and misclassification of the current regime.
    From an applied perspective, this framework can model the monitoring of radar signals with different disturbance levels. For example, $\theta_t = 0$ may indicate that the system is ``off,'' while $\theta_t = k$ may indicate a level-$k$ disruption. The system switches between ``on'' and ``off,'' with different disturbance levels in the ``on'' state. The observer seeks to identify signal changes as accurately and quickly as possible. For example, the penalty function may take the form
    \[c(d,\theta) = 1_{d=\theta} c^1(\theta) +1_{d\neq\theta}c^2(\theta),\]
    where $c^2\geq c^1$, so an incorrect decision incurs a larger penalty. By Theorem \ref{thm_suff_1}, the generator is hypoelliptic. This model can be viewed as the Brownian counterpart of the hidden-regime problem studied in detail in \cite{BL}.
\end{example}
\begin{example}(Byzantine testing and detection)\\
\label{eg4}
We consider a quickest detection problem with channels that may be initially faulty. This is a Bayesian version of the formulation in \cite{BL_byzantine}. Suppose that we observe two processes $X^1$ and $X^2$, driven by independent one-dimensional Brownian motions $W^1$ and $W^2$, respectively. Assume that there is a disorder time $\eta$, independent of $W^1$ and $W^2$, with distribution
\[\mathbb{P}_0(\eta=0)=p,\ \ \ \mathbb{P}_0(\eta>t)=(1-p)e^{-\lambda t},\ \ 
t\geq 0\]
for some $p\in[0,1]$. The drift of $X^1$ takes values in $\{\mu_0,\mu_1\}$, while the drift of $X^2$ takes values in $\{m_0,m_1\}$. At time $\eta$, process $X^1$ changes its drift from $\mu_0$ to $\mu_1$, and process $X^2$ changes its drift from $m_0$ to $m_1$, if either process has not already done so. Our goal is to detect the change time $\eta$. \par
However, the initial model may itself be misspecified: the two channels may not have zero initial drifts. We call this example ``Byzantine'' because it is related to the Byzantine generals problem. Suppose that the channels may be corrupted at time $0$ and behave as though they have already acquired positive drifts. An observation may therefore come from an untrustworthy channel, increasing the risk of a false alarm if it is used indiscriminately. Specifically, let $\xi$ denote the initial corruption state, independent of $W^1$ and $W^2$. We assume that $\xi$ can take four values represented by the initial drift configuration:
\begin{equation}
\xi^T = 
\begin{cases}
& (\mu_1,m_1),\ \ \ \text{both channels are affected at $t=0$}\\
& (\mu_1,m_0),\ \ \ \text{$X^1$ is affected at time $t=0$}\\
& (\mu_0,m_1),\ \ \ \text{$X^2$ is affected at time $t=0$}\\
& (\mu_0,m_0),\ \ \ \text{no channels are affected at time $t=0$},
\end{cases}
\end{equation}
where $\mu_1\neq \mu_0$ and $m_1\neq m_0$. 
We assign the prior probabilities $\pi=(\pi_0, \pi_1,\pi_2,\pi_3)$ to these four states, respectively. 
Let $\theta_t$ denote the current drift configuration. Before $\eta$ it is determined by $\xi$, while at time $\eta$ any channel that has not already acquired its affected drift does so. Thus, after $\eta$, the drift
configuration is $(\mu_1,m_1)$. We consider the detection problem
\[V = \inf_{\tau \in \mathcal{T}}\mathbb{E}[(\eta-\tau)^+]+c\mathbb{E}[(\tau-\eta)^+].\]
Let $P_t$ denote the posterior probability that the disorder has not yet occurred. Specifically, $P_t=\P(\eta>t\vert \F_t)=\sum_{i=0}^3\tilde\Pi_t^i$, where $\tilde\Pi_t^i:=\P(\xi=i,\eta>t\vert\F_t)$. The problem has two notable features. First, the drifts of the processes are highly correlated through the disorder time $\eta$. Second, we must infer whether the channels are corrupted and therefore how much to trust each one. Conditional on the state $\theta$, the drift-difference matrix $A$ has the form
\[A= -\begin{bmatrix}
    \mu_1-\mu_0& 0&\mu_1-\mu_0\\
    0&m_1-m_0&m_1-m_0
\end{bmatrix}.\]
By the memoryless property of the exponential disorder time, the current drift configuration $\theta$ has generator
\[Q= \begin{bmatrix}
     0&0&0&0\\
    \lambda&-\lambda&0&0\\
    \lambda&0&-\lambda&0\\
    \lambda&0&0&-\lambda
\end{bmatrix}.\]
The matrix $A$ has rank $2$, but the vector $(\|a_1\|^2, \|a_2\|^2, \|a_3\|^2)$ always belongs to its row space. However, for $j = 1,2,3$, $q_{mj} = 0$ for all $m\neq j$, while $q_{10} = q_{20} = q_{30} = \lambda>0$, so the assumptions of Proposition \ref{thm_suff_2} are satisfied. The vector $(1,1,1)$ does not belong to the row space of $A$ because $\mu_1\neq \mu_0$ and $m_1\neq m_0$. Consequently, Proposition \ref{thm_suff_2} shows that the posterior operator associated with the four-state drift-configuration filter satisfies H\"ormander's condition on $(0,\infty)^3$. Since $\rank[A^T\vert 1]=3$, the parabolic H\"ormander condition also holds. If we decide that state $1$, $2$, or $3$ is true, we may choose to disregard information from certain coordinates, so the cost function may depend on both $\pi$ and $x$. However, since $n=3$, the joint operator $\calL_{\phi,x}$ fails to satisfy H\"ormander's condition. This example shows that the parabolic H\"ormander condition does not imply the joint H\"ormander condition: augmenting the state with time and augmenting it with the observation process are fundamentally different operations. This problem has many potential applications, and its detailed properties are of independent interest. \end{example}

\bibliographystyle{abbrv}
\nocite*
\bibliography{references} 

\end{document}